\documentstyle[12pt]{article}
\makeatletter

\@addtoreset{equation}{section}

\def\section#1{\@ifstar{Nsection*}{\def\thesection{\arabic{section}.}\@startsection{section}{1}{\z@}{-3.5ex
plus -1ex minus
 -.2ex}{2.3ex plus .2ex}{\bf}{\uppercase\expandafter{#1}}
\def\thesection{\arabic{section}}}}

\def\subsection{\@startsection{subsection}{2}{\z@}{-3.25ex plus -1ex minus
 -.2ex}{-1em}{\bf}}

\def\thebibliography#1{\par\vskip 6 mm\noindent {\bf REFERENCES} \par \list
 {[\arabic{enumi}]}{\settowidth\labelwidth{[#1]}\leftmargin\labelwidth
 \advance\leftmargin\labelsep
 \usecounter{enumi}}
 \def\newblock{\hskip .11em plus .33em minus -.07em}
 \sloppy
 \sfcode`\.=1000\relax}

\tabbingsep \labelsep

\def\address#1{\gdef\@address{#1}}

\def\maketitle{\par
 \begingroup
 \def\thefootnote{\fnsymbol{footnote}}
 \def\@makefnmark{\hbox
 to 0pt{$^{\@thefnmark}$\hss}}
 \if@twocolumn
 \twocolumn[\@maketitle]
 \else \newpage
 \global\@topnum\z@ \@maketitle \fi\thispagestyle{empty}\@thanks
 \endgroup
 \setcounter{footnote}{0}
 \let\maketitle\relax
 \let\@maketitle\relax
 \gdef\@thanks{}\gdef\@author{}\gdef\@title{}\let\thanks\relax}
\def\@maketitle{\newpage
 \null
 \vskip 10mm \begin{center}
 {\bf\uppercase\expandafter{\@title}\par} \vskip 4 mm 
 {\bf\@author \par}
 \end{center}
 \vskip 20 mm}

\newcommand{\nc}{\newcommand}
\newtheorem{lemma}{Lemma} [section]
\newtheorem{theorem}[lemma]{Theorem}
\newtheorem{proposition}[lemma]{Proposition}
\newtheorem{corollary}[lemma]{Corollary}
\newtheorem{remark}[lemma]{Remark}
\newtheorem{assumption}[lemma]{Assumption}
\newtheorem{example}[lemma]{Example}

\newtheorem{definition}[lemma]{Definition}

\nc{\QED}{\mbox{}\hfill \raisebox{-2pt}{\rule{5.6pt}{8pt}\rule{4pt}{0pt}} 
          \medskip\par}

\nc{\di}{\displaystyle}
\nc{\SS}{{\rm I\mkern-4.0mu S}}
\nc{\SSs}{{\sf S\mkern-6.5mu S}} 
\nc{\C}{{\rm I\mkern-4.0mu C}}        
\nc{\R}{{\rm I\mkern-4.0mu R}}
\nc{\Z}{{\sf Z\mkern-6.5mu Z}}
\nc{\CCa}{\mbox{$C \mkern-10.5mu\raisebox{.18em}{$\scriptstyle/$}\mkern3mu$}}
\nc{\CCb}{C\mkern-11.5mu I\mkern4mu}
\nc{\Cc}{\mathchoice{\CCa}{\CCa}{\CCb}{C\mkern-10.5mu I\mkern3mu}}

\nc{\Ref}[1]{~\mbox{\rm (\ref{#1})}}       
\nc{\REF}[1]{~\mbox{\rm \ref{#1}}}             
\nc{\Norm}[2]{\|#1\|\left.\vphantom{T_{j_0}^0}\!\!\right._{#2}}                    
                                               
\nc{\Log}{\mathop {\rm Log}\nolimits}

\expandafter\chardef\csname pre amssym.def at\endcsname=\the\catcode`\@
\catcode`\@=11

\def\undefine#1{\let#1\undefined}
\def\newsymbol#1#2#3#4#5{\let\next@\relax
 \ifnum#2=\@ne\let\next@\msafam@\else
 \ifnum#2=\tw@\let\next@\msbfam@\fi\fi
 \mathchardef#1="#3\next@#4#5}
\def\mathhexbox@#1#2#3{\relax
 \ifmmode\mathpalette{}{\m@th\mathchar"#1#2#3}%
 \else\leavevmode\hbox{$\m@th\mathchar"#1#2#3$}\fi}
\def\hexnumber@#1{\ifcase#1 0\or 1\or 2\or 3\or 4\or 5\or 6\or 7\or 8\or
 9\or A\or B\or C\or D\or E\or F\fi}

\font\tenmsa=msam10
\font\sevenmsa=msam7
\font\fivemsa=msam5
\newfam\msafam
\textfont\msafam=\tenmsa
\scriptfont\msafam=\sevenmsa
\scriptscriptfont\msafam=\fivemsa
\edef\msafam@{\hexnumber@\msafam}
\mathchardef\dabar@"0\msafam@39
\def\dashrightarrow{\mathrel{\dabar@\dabar@\mathchar"0\msafam@4B}}
\def\dashleftarrow{\mathrel{\mathchar"0\msafam@4C\dabar@\dabar@}}

\def\ulcorner{\delimiter"4\msafam@70\msafam@70 }
\def\urcorner{\delimiter"5\msafam@71\msafam@71 }
\def\llcorner{\delimiter"4\msafam@78\msafam@78 }
\def\lrcorner{\delimiter"5\msafam@79\msafam@79 }
\def\yen{{\mathhexbox@\msafam@55 }}
\def\checkmark{{\mathhexbox@\msafam@58 }}
\def\circledR{{\mathhexbox@\msafam@72 }}
\def\maltese{{\mathhexbox@\msafam@7A }}

\font\tenmsb=msbm10
\font\sevenmsb=msbm7
\font\fivemsb=msbm5
\newfam\msbfam
\textfont\msbfam=\tenmsb
\scriptfont\msbfam=\sevenmsb
\scriptscriptfont\msbfam=\fivemsb
\edef\msbfam@{\hexnumber@\msbfam}
\def\Bbb#1{\fam\msbfam\relax#1}
\def\widehat#1{\setbox\z@\hbox{$\m@th#1$}%
 \ifdim\wd\z@>\tw@ em\mathaccent"0\msbfam@5B{#1}%
 \else\mathaccent"0362{#1}\fi}
\def\widetilde#1{\setbox\z@\hbox{$\m@th#1$}%
 \ifdim\wd\z@>\tw@ em\mathaccent"0\msbfam@5D{#1}%
 \else\mathaccent"0365{#1}\fi}
\font\teneufm=eufm10
\font\seveneufm=eufm7
\font\fiveeufm=eufm5
\newfam\eufmfam
\textfont\eufmfam=\teneufm
\scriptfont\eufmfam=\seveneufm
\scriptscriptfont\eufmfam=\fiveeufm

\catcode`\@=\csname pre amssym.def at\endcsname

\begin{document}

\title{A class of pseudo-differential operators with oscillating symbols} 
\author{ D.Yafaev
\\Universit\'e de Rennes}
\maketitle

\begin{abstract}
We study a class of  pseudo-differential operators  with oscillating symbols or oscillating 
amplitudes appearing in the long-range scattering theory. We develop the basic calculus for
operators  from such classes and solve some concrete problems posed by applications to  
scattering theory, especially to the scattering matrix. In particular, we show  that under natural
assumptions the spectrum of a pseudo-differential operator with an oscillating symbol covers the
unit circle.
\end{abstract}

\null\hspace{7mm}{\bf Contents}
\newline\null\hspace{2mm} 1. Introduction
\newline\null\hspace{2mm} 2. The basic calculus
\newline\null\hspace{2mm} 3. The action on an exponent
\newline\null\hspace{2mm} 4. The essential spectrum
\newline\null\hspace{2mm} 5. Integral kernels

\section {Introduction}

Pseudo-differential operators  (PDO)  $A$,
\begin{equation}
  (A u)(x)=(2\pi)^{-d/2}\int_{{\Bbb R}^d}e^{i <x,\xi> } a(x,\xi)
\hat{u}(\xi) d\xi,
\label{eq:BC1}\end{equation}
or, more generally,
\begin{equation}
 (Au)(x)=(2\pi)^{-d}\int_{X} \int_{{\Bbb R}^d} e^{i<x-x^\prime,\xi>} {\bf a}(x,x^\prime,\xi)
u(x^\prime) dx^\prime d\xi
\label{eq:RH1}\end{equation}
are well-defined as mappings 
$A:C_0^\infty(X)\rightarrow C^\infty(X)$
for symbols $a$ or amplitudes ${\bf a}$ from the H\"ormander classes ${\cal S}^m_{\rho,\delta}$ or
${\cal S}^m_{\rho,\delta,\delta}$ for arbitrary $\rho>0$ and $\delta<1$. Here $X\subset {\Bbb R}^d$
is an open set,  $\hat{u}$ is the Fourier transform of a function $u\in C_0^\infty(X)$ and the
definition of the   classes ${\cal S}^m_{\rho,\delta}$ and
${\cal S}^m_{\rho,\delta,\delta}$  is recalled in subsection 2.1. A crucial advantage of the PDO
theory is that rather an advanced calculus can be developed  (see \cite{Hor3},
\cite{Sh}) in its framework. For example, one obtains formulas for symbols of the  adjoint operator
$A^\ast$, of the product $A_1 A_2$ of two PDO, finds a relation between an amplitude and the
corresponding symbol, checks the invariance of the theory with respect to  change of variables and
so on. Such an advanced calculus can be conveniently developed in the 
  classes ${\cal S}^m_{\rho,\delta}$ and
${\cal S}^m_{\rho,\delta,\delta}$  which, however, requires the assumption $\rho>1/2>\delta$.

Our aim is to study a class of  pseudo-differential operators  with oscillating symbols or
oscillating  amplitudes appearing  (see \cite{Y}) in the long-range scattering theory. More precisely, we consider
PDO (\ref{eq:BC1}) and (\ref{eq:RH1}), where
\begin{equation}
 a(x,\xi)=e^{i \Phi  (x,\xi) } b(x,\xi), \quad  \Phi\in {\cal S}^r,\quad b\in {\cal S}^m,\quad
r\in[0,1),
\label{eq:BC2}\end{equation}
or
\begin{equation} 
{\bf a}(x,x^\prime,\xi)=e^{i \Theta  (x,x^\prime,\xi) } {\bf b}(x,x^\prime,\xi),
\quad  \Theta\in {\cal S}^r,\quad {\bf b}\in {\cal S}^m,\quad r\in[0,1),
\label{eq:I1}\end{equation}
and $  {\cal S}^m= {\cal S}^m_{1,0}$ or $  {\cal S}^m= {\cal S}^m_{1,0,0}$. Let us denote the
classes of symbols (\ref{eq:BC2}) and amplitudes (\ref{eq:I1}) by ${\cal C}^m(\Phi)$ and ${\cal
C}^m(\Theta)$, respectively. Sometimes we use the same notation  ${\cal C}^m(\Phi)$, ${\cal
C}^m(\Theta)$, ${\cal S}^m_{\rho,\delta}$ and
${\cal S}^m_{\rho,\delta,\delta}$  for PDO (\ref{eq:BC1}) and (\ref{eq:RH1}) with symbols and
amplitudes from the corresponding classes. Clearly,
\[
{\cal C}^m(\Phi)\subset {\cal S}^m_{1-r,r},\quad {\cal C}^m(\Theta)\subset {\cal S}^m_{1-r,r,r}
\]
so that (\ref{eq:BC2}) and  (\ref{eq:I1}) are ``good" classes if
$r<1/2$. On the other hand, the standard calculus fails for operators from these classes if $r\geq
1/2$. In this paper we consider several  concrete problems 
for PDO from classes ${\cal C}^m(\Phi)$ and ${\cal C}^m(\Theta)$ posed
by the long-range scattering theory.

\medskip

Of course, standard formulas of the PDO calculus  for the adjoint
$A^\ast$ or for the product $A_1 A_2$ fail in the class ${\cal C}^m(\Phi)$
if $\Phi\in {\cal S}^r$ with $r\in[1/2,1)$. Fortunately, in applications to   scattering theory
only the combinations $A_1 A_2^\ast$ and $ A_2^\ast A_1$ appear. In  Section 2 we show that
\begin{equation}
 {\rm if}\;  A_j\in {\cal C}^{m_j}(\Phi),\; j=1,2,\quad {\rm then}\; A_1 A_2^\ast\in{\cal S}^m,\; 
 A_2^\ast A_1\in{\cal S}^m, \;  {\rm where}\; m=m_1+m_2.
\label{eq:I1a}\end{equation}
We justify also usual expansions for symbols of the operators $A_1 A_2^\ast$ and $ A_2^\ast A_1$;
in particular, their  principal symbols are equal $b_1\overline{b_2}$.
 Note that inclusions (\ref{eq:I1a})  were checked by a different method in
\cite{K} but, to best of our knowledge, asymptotic expansions  for   symbols of $A_1 A_2^\ast$ and
 $A_2^\ast A_1$ are new.

Each of the inclusions (\ref{eq:I1a}) imply that a PDO $A$ with  compactly supported symbol
$a\in{\cal C}^0(\Phi)$ is bounded in the space $L_2(X)$. In the case $r\leq 1/2$ this  follows
 from  results of \cite{Ha} and
\cite{CV} for PDO from classes ${\cal S}^0_{\rho,\delta}$ with $\rho \geq \delta$ but, if $r
> 1/2$, then these general results can no longer be applied. 

 We show  also that a PDO $A$ defined by its amplitude ${\bf a}\in {\cal C}^m(\Theta)$ admits
representation  (\ref{eq:BC1}) and find
an expression for the symbol $a$ of this PDO in terms of the amplitude ${\bf a}$.
It is different from the familiar expression in the case 
${\bf a}\in {\cal S}^m_{\rho,\delta,\delta}$ 
 with $\rho>\delta$; in particular, the ``principal" symbol of $A$ does not coincide with ${\bf
a}(x,x,\xi)$. On the other hand, if $\Theta(x,x,\xi)=0$, then  PDO (\ref{eq:RH1}) with  
amplitude (\ref{eq:I1}) has   symbol $a\in{\cal S}^m$; moreover, $a$ admits the usual expansion in
terms of the amplitude ${\bf a}$.

\medskip

In Section 3 we calculate the action of a PDO $A$ on an exponent
$u_{\lambda}(x)=e^{i\lambda\psi(x)} f(x)$, where $f\in C_0^\infty (X)$ and
$\lambda\rightarrow\infty$. This result is used in the following section for construction of Weyl
(singular) sequences. In view of this application to the spectral theory, we have to
study the case when $f$ depends on an additional parameter $\varepsilon$ 
 and supports of functions $f_\varepsilon$ are shrinking to
some point $x_0\in X$; moreover, the function $\psi$ may depend  on $\lambda$. It suffices for us to
find only the leading term of $A u_{\lambda,\varepsilon}$ and  give an estimate of the rest in the space $L_2$.
Our calculation goes through if the localization in $x$ is not too sharp compared to $\lambda^{-1}$.
 More precisely, we show in Section 3 that if
\begin{equation}
 u_{\lambda,\varepsilon}(x)=e^{i\lambda\psi(x,\lambda)} f_\varepsilon(x)
\label{eq:I2}\end{equation}
and $ f_\varepsilon(x)=\varepsilon^{-d/2}f((x-x_0)/\varepsilon)$,
then 
\begin{equation}
 (A u_{\lambda,\varepsilon})(x)=e^{iG(x,\lambda)} b(x ,
\lambda\psi^\prime(x,\lambda)) u_{\lambda,\varepsilon}(x) + o(1),
\label{eq:I3}\end{equation}
 where the function $G(x,\lambda)$ is determined by $\Phi$ and
$\psi$  and $o(1)$ tends to zero in $L_2$ as $\lambda\rightarrow\infty$ and
$\lambda^{1-r}\varepsilon\rightarrow\infty$.

\medskip
 
In Section 4 we study the essential  spectrum of a PDO $A\in{\cal C}^0(\Phi)$ in the space $L_2(X)$.
A typical result is the following. Suppose that for some point $x_0\in X,\; \xi_0\neq 0$
\begin{equation}
\lim_{\lambda\rightarrow\infty} b(x_0, \lambda\xi_0)=\mu_0 \neq 0.
\label{eq:I1b}\end{equation}
Then, under some mild assumptions on the phase function
$\Phi$, the spectrum of the operator (\ref{eq:BC1}) with symbol (\ref{eq:BC2}) covers the whole 
circle ${\Bbb T}_{\kappa}=\{z\in{\Bbb C}:|z|=\kappa\}$, where $\kappa =|\mu_0|$. In particular,  the
spectrum of
$A$ covers the unit  circle if, for example, $b=1$.

For the proof, we construct, for any point $\mu\in{\Bbb T}_\kappa$, a Weyl sequence which
we seek in the form (\ref{eq:I2}) where 
 $\lambda\rightarrow\infty,\varepsilon\rightarrow 0$ but $\lambda^{1-r}\varepsilon\rightarrow
\infty$. We  construct $\psi$ in such a way that $G(x,\lambda)$
essentially does not depend on $x$ in a neighbourhood of the point $x_0$ so that $G(x,\lambda)$ may
be replaced by $G(x_0,\lambda)$ in (\ref{eq:I3}). If $r<1/2$  (this case was considered in \cite{Y}),
we can set 
\begin{equation}
\psi(x)=<\xi_0, x-x_0>
\label{eq:1.9}\end{equation}
but in the general case $\psi$ is a polynomial of degree $[r(1-r)^{-1}]+1$
 with coefficients depending on $\lambda$.
 If $|G(x_0,\lambda)|\rightarrow\infty$ as
$\lambda\rightarrow\infty$, then it is possible,  for any  $\mu_1\in{\Bbb T}_1$, to find  a 
sequence $\lambda_p\rightarrow\infty$ such that $ e^{iG(x_0,\lambda_p)} =\mu_1$.
Then, for a suitable sequence $\varepsilon_p\rightarrow 0$, $u_{\lambda_p,\varepsilon_p}$ is a Weyl
sequence for the operator $A$ and the point $\mu=\mu_1\mu_0$.

\medskip

In Section 5 we consider PDO as integral operators in direct integrals of multiplication operators.
For example, passing to the spherical coordinates and denoting $\lambda=x^2$ we can represent
$L_2({\Bbb R}^d)$ as the space
$L_2({\Bbb R}_+; L_2({\Bbb S}^{d-1}))$ of vector-functions. This gives the spectral representation of
the multiplication by $x^2$. In this representation an operator $A$ can be considered as a formal
integral operator, whose kernel $A^\natural (\mu,\nu)$ is, for every $\mu,\nu>0$, an operator in
the space $L_2({\Bbb S}^{d-1})$. A precise definition of the kernel requires,
of course, some assumptions on the operator $A$. Suppose, for example, that, for some $s>1/2$, an
operator $A$ is bounded from the Sobolev space $H^{-s}({\Bbb R}^d)$ into the space $H^s({\Bbb
R}^d)$. In this case $A^\natural (\mu,\nu)$  is well-defined as a bounded operator in the space
$L_2({\Bbb S}^{d-1})$ and is a continuous function of $\mu,\nu>0$. This implies that the same
result holds for a PDO $A$ if its symbol $a$ belongs to the class ${\cal S}^m_{1,0}$ with $m<-1$.
In this case $A^\natural (\mu,\nu)$ is also a PDO and one can give (see \cite{BY}) an explicit
expression for its principal symbol. The case $m\geq -1$ was studied in \cite{LY}. In particular,
for $m=-1$ it was shown there that the diagonal value $A^\natural (\lambda,\lambda)$  is correctly
defined if (and only if) the principal symbol $a_{-1}(x,\xi)$ of $A$ equals to zero on the conormal
bundle to the sphere $\{|x|=\lambda^{1/2}\}$, that is $a_{-1}(\lambda^{1/2} \omega,t\omega)=0$ for
$|\omega|=1$ and $t\in{\Bbb R}$ (for sufficiently large $|t|$). 

Our goal in Section 5 is to consider PDO with symbols from  arbitrary classes ${\cal
S}^m_{\rho,\delta}$. In such a general case there is no invariance with respect to  change of
variables, and hence we are obliged to work with PDO defined by their amplitudes ${\bf a}\in{\cal
S}^m_{\rho,\delta,\delta}$ where    $\rho>0$ and $\delta<1$ are arbitrary.
So,  compared to \cite{LY}, we consider  PDO from a more general class but our
condition  on the conormal bundle is much more restrictive. Actually, we suppose that
 ${\bf a}(x,x^\prime,\xi)=0$ for $x^2$ close to some $\lambda_0>0$, small $|x-x^\prime|$ and $\xi$
from some conical neighbourhood of the  line $tx$, $t\in{\Bbb R}$. Then we 
 construct a continuous
kernel $A^\natural (\mu,\nu)$ in a neighbourhood of the point $\mu=\nu=\lambda_0$. We also check that
$A^\natural (\mu,\nu)$ is  a PDO from the class ${\cal S}^{m+1}_{\rho,\delta,\delta}$ and give
an explicit expression for its amplitude.

 Our result on the existence of  diagonal values 
$A^\natural (\lambda,\lambda)$ is, to a certain extent, similar in spirit to a result of
\cite{Hor3}, Chapter 8, on the existence of restriction of a distribution  to a manifold ${\bf S}$.
In \cite{Hor3} it is required that the wave front of the distribution does not intersect with the
conormal bundle to ${\bf S}$. In the example above, this implies that kernel $k(x,x^\prime)$ of a
PDO $A$ can, in some sense, be restricted to
${\bf S}=S\times S$ where $S$ is the sphere  $|x|=\lambda^{1/2}$.
\medskip

 The results of  Section 2 are a necessary technical background for an
elementary proof of the existence and completeness of wave operators for the Schr\"odinger operator
with a long-range potential. The results of the following sections are used for a study of
the corresponding scattering matrix. Actually, the singular part of the scattering matrix is
defined as a diagonal value of kernel of some PDO. This requires the results of  Section 5. A
study of spectral properties of the scattering matrix relies on the results of Sections 3 and 4.

\section { the basic calculus}

{\bf 2.1.} We recall first the definition of the H\"ormander classes ${\cal S}^{m}_{\rho,\delta}$.
 Let $X\subset  {\Bbb R}^{d}$ be some open set and let $m\in {\Bbb R}$. The set ${\cal
S}^{m}_{\rho,\delta}={\cal S}^{m}_{\rho,\delta}(X\times {\Bbb R}^d)$ consists of functions
 $a \in C^\infty (X\times {\Bbb R}^d)$ such that, for all multi-indices $\alpha,\beta$ and all
compact
$K\subset X$, there exist $ N_{\alpha,\beta,K}$ such that
\[ 
 |(\partial_{\xi}^{\alpha} \partial_{x}^{\beta} a) (x, \xi)  |\leq N_{\alpha,\beta,K}(a) (1
+|\xi|)^{m-|\alpha|\rho+ |\beta| \delta}
\]
 for all $(x, \xi) \in K\times {\Bbb R}^{d}$.
  The best $N_{\alpha,\beta,K}(a)$ are the semi-norms of the symbol $a$.  The set ${\cal
S}^{m}_{\rho,\delta,\delta}$ of functions ${\bf a}(x,x^\prime,\xi)$ is defined exactly in the same
way if $x$ is replaced by $(x,x^\prime)$.  We denote
${\cal S}^{m}={\cal S}^{m}_{1,0}$ or ${\cal S}^{m}={\cal S}^{m}_{1,0,0}$.
 Below $C$ and $c$ are different  positive constants whose values are unimportant.

For any PDO (\ref{eq:RH1}) with the amplitude  ${\bf a}\in {\cal S}^{m}_{\rho,\delta,\delta}$, where
$\rho>0,\delta<1$, its kernel
\begin{equation}
 k(x,x^\prime)=k_A(x,x^\prime)=(2\pi)^{-d}  \int_{{\Bbb R}^d} e^{i<x-x^\prime,\xi>} {\bf
a}(x,x^\prime,\xi)  d\xi
\label{eq:RH1b}\end{equation}
 is well-defined for $x\neq x^\prime$ and $k(x,x^\prime)$
 is a $C^\infty$-function outside of the diagonal
$x=x^\prime$. A PDO  is called properly supported if $k(x,x^\prime)=0$ for any $x$ and $|x^\prime|$
 sufficiently large, i.e. $|x^\prime|\geq c(x)$, as well as for any $x^\prime$ and  $|x|\geq
c(x^\prime)$. Any properly supported PDO (\ref{eq:RH1}) can be written in the form (\ref{eq:BC1})
with  symbol
\begin{equation}
  a(x,\xi)= (2\pi)^{-d} \int_{{\Bbb R}^d} \int_{{\Bbb R}^d} {\bf a}(x,x+z,\xi+\zeta) e^{-i <z,\zeta>
} dzd\zeta
\label{eq:BC7}\end{equation}
 but, of course, in general it cannot be claimed that $a\in{\cal
S}^{m}_{\rho,\delta}$ if ${\bf a}\in{\cal S}^{m}_{\rho,\delta,\delta}$.

 A standard integration by parts in the variable $z$ based on the formula
\begin{equation}
 e^{-i <z,\zeta>} = \langle \zeta\rangle^{-k} \langleÊD_z\rangle ^k e^{-i <z,\zeta>}, 
\quad \langleÊD_z\rangle ^2= I+ D^2_{z_1}+\ldots + D^2_{z_d}, \quad k \; {\rm is}\; {\rm even},
\label{eq:IP}\end{equation}
 shows only that (for any  $\rho\geq 0,\delta <1$)
\begin{equation}
 |\partial_\xi^\alpha \partial_x^\beta a(x,\xi)| \leq C (1+|\xi|)^n,
\quad (x, \xi) \in K\times {\Bbb R}^{d},
\label{eq:BC7a}\end{equation}
 for all multi-indices $\alpha,\beta$, all compact
$K\subset X$ and some  $n=n(\alpha,\beta)$,
$C=C( \alpha,\beta,K)$.

{\bf 2.2.} 
The following result shows that the class ${\cal C}^m (\Theta)$ reduces to the ``best"
class ${\cal S}^m $ if the phase function $\Theta(x,x^\prime,\xi)$ equals to zero at the diagonal
$x=x^\prime$.

\begin{theorem}\label{AM} 
 Suppose that a PDO $A$ is given by formula $(\ref{eq:RH1})$ where ${\bf a}$ admits representation
$(\ref{eq:I1})$ and $\Theta(x,x ,\xi)=0$. Then
$A$ is a PDO with symbol $a\in {\cal S}^m$  and 
$a(x,\xi)$ admits the asymptotic expansion
\begin{equation}
  a(x,\xi)=\sum_{|\alpha|\geq 0} (\alpha !)^{-1} a_\alpha (x,\xi),
\quad {\rm where}\quad
  a_\alpha(x,\xi)=(\partial^\alpha_\xi D^\alpha_{x^\prime} {\bf a})(x,x^\prime,\xi)|_{x=x^\prime};
\label{eq:AM1}\end{equation}
 in particular, $a_\alpha\in{\cal S}^{m-|\alpha|(1-r)}$ for all $\alpha$.
\end{theorem}

Remark first of all that 
\begin{equation}
 (D^\alpha_{x^\prime} {\bf a})(x,x^\prime,\xi)= e^{i \Theta  (x,x^\prime,\xi) } {\bf b}_\alpha
(x,x^\prime,\xi),\quad {\bf b}_\alpha \in  {\cal S}^{m +|\alpha|r},
\label{eq:BC2a}\end{equation}
 so that the oscillating factor disappears if $x=x^\prime$ and, consequently, $a_\alpha\in{\cal
S}^{m-|\alpha|(1-r)}$. Thus the series (\ref{eq:AM1}) is asymptotic for any $r<1$.

Expansion (\ref{eq:AM1}) for  the symbol $a$ in terms of the amplitude ${\bf a}$ is, of course, the 
 same  as in the case ${\bf a} \in {\cal S}^m_{\rho,\delta,\delta}$ for $\rho>\delta$. Thus,
Theorem~\ref{AM} holds true in the case $\Theta\in{\cal S}^r$ for $r<1/2$. Since the right-hand side
of (\ref{eq:AM1}) is an asymptotic series for any $r<1$, this makes quite plausible that equality
(\ref{eq:AM1}) itself is also valid for all $r<1$. Some part of the usual proof (expounded, for
example, in \cite{Sh}) of (\ref{eq:AM1}) for 
${\bf a}\in {\cal S}^m_{\rho,\delta,\delta} $ with $\rho>\delta$ applies for arbitrary $\rho>0$,
$\delta<1$. We try to avoid repeating these arguments and concentrate on the part of the construction
which fails for $\rho\leq \delta$.

At a formal level (\ref{eq:AM1})  is a consequence of the following elementary

\begin{lemma}\label{Ta}
Let ${\bf p}(z,\zeta)$ be a $C^\infty$-function which is compactly supported
in the variable
$z$ and is polynomially bounded in the variable $\zeta$. Then for any $N\geq 1$
\begin{equation}
 (2\pi)^{-d} \int_{{\Bbb R}^d} \int_{{\Bbb R}^d} {\bf p}(z,\zeta) e^{-i <z,\zeta> }dzd\zeta  =
\sum_{0\leq |\alpha|\leq N-1} (\alpha!)^{-1}(\partial^\alpha_\zeta D^\alpha_z{\bf p}) 
(0,0)  + p^{(N)},
\label{eq:Ta2}\end{equation}
  where
\begin{equation}
 p^{(N)} =(2\pi)^{-d} N \sum_{|\alpha| =N} (\alpha!)^{-1} 
\int_0^1 (1-t)^{N-1}
\Bigl(\int_{{\Bbb R}^d} \int_{{\Bbb R}^d} (\partial^\alpha_\zeta {\bf p})(z,t\zeta)
\zeta^\alpha e^{-i <z,\zeta> } dzd\zeta\Bigr) dt. 
\label{eq:Ta3}\end{equation}
\end{lemma}
 {\it Proof.} --
 Let us use the Taylor expansion (with the rest) at the point $\zeta=0$ 
\[
 {\bf p}(z,\zeta) =\sum_{|\alpha|\leq N-1} (\alpha!)^{-1}(\partial^\alpha_\zeta {\bf p}) 
(z,0)\zeta^\alpha + {\bf\tilde{p}}^{(N)} (z,\zeta),
\]
 where
\[
 {\bf\tilde{p}}^{(N)} (z,\zeta) = N \sum_{|\alpha| =N} (\alpha!)^{-1}\zeta^\alpha 
\int_0^1 (1-t)^{N-1}(\partial^\alpha_\zeta {\bf p})  (z,t\zeta)dt. 
\]
 Taking into account that
\[
  (2\pi)^{-d}\int_{{\Bbb R}^d} \int_{{\Bbb R}^d} (\partial^\alpha_\zeta {\bf p}) 
(z,0) 
 e^{-i <z,\zeta>} \zeta^\alpha dz d\zeta = ( D^\alpha_z \partial^\alpha_\zeta {\bf p})  (0,0)
\]
 we arrive at equality (\ref{eq:Ta2}).$\quad\Box$

Our proof of Theorem~\ref{AM} relies, of course, on  representation (\ref{eq:BC7}) for the symbol
$a(x,\xi)$. Let us use, for any fixed $x,\xi$,  Lemma~\ref{Ta} with 
\[
 {\bf p}(z,\zeta;x,\xi) = {\bf a}(x,x+z,\xi+\zeta).
\]
 Then the coefficients $( D^\alpha_z \partial^\alpha_\zeta {\bf p}) 
(0,0;x,\xi)$ coincide with the numbers
$a_\alpha(x,\xi)$ defined in (\ref{eq:AM1}). Equality (\ref{eq:Ta2}) shows that $a(x,\xi)$ is a sum
of  the terms $(\alpha!)^{-1} a_\alpha(x,\xi)$ over $\alpha,\; 0\leq |\alpha|\leq N-1$, and of the
rest $a^{(N)}(x,\xi)$. Integrating by parts in (\ref{eq:Ta3}) we obtain that
\[
a^{(N)}(x,\xi) =(2\pi)^{-d} N \sum_{|\alpha| =N} (\alpha!)^{-1} 
\int_0^1 (1-t)^{N-1} R^{(\alpha)} (x,\xi;t) dt, 
\]
where 
\begin{equation}
 R^{(\alpha)} (x,\xi;t) =
\int_{{\Bbb R}^d} \int_{{\Bbb R}^d}(\partial^\alpha_\xi D^\alpha_{x^\prime}{\bf a})  (x,x+z,\xi
+t\zeta)  e^{-i <z,\zeta>} dz d\zeta.
\label{eq:BC10}\end{equation}
 Taking into account (\ref{eq:BC7a}) we see that for the proof of Theorem~\ref{AM}  it suffices to
check the following

\begin{lemma}\label{AM1}
 Under the assumptions of Theorem~$\ref{AM}$ there exists a number $q$ such that
 for all $\alpha$ and all compact $K\subset X$
\begin{equation}
 |R^{(\alpha)} (x,\xi;t)| \leq C  (1+|\xi|)^{q-|\alpha|(1-r)},
\quad   (x, \xi) \in K\times {\Bbb R}^{d},
\label{eq:AM3}\end{equation} 
 uniformly in $t\in[0,1]$.
\end{lemma}

An estimation of integral (\ref{eq:BC10}) is different for large and small values of $|\zeta|\,
|\xi|^{-1}$. In the first case it is quite standard. The following lemma is true for an arbitrary
amplitude
${\bf a}\in {\cal S}^m_{\rho,\delta,\delta}$ with any $\rho\geq 0$, $\delta<1$. Its proof relies on
equality (\ref{eq:IP}) and integration by parts in the variable $z$. 

\begin{lemma}\label{AM11}
If $|\xi|\leq C |\zeta|$, then for any $n$
\[
|\int_{{\Bbb R}^d}(\partial^\alpha_\xi D^\alpha_{x^\prime}{\bf a})  (x,x+z,\xi
+t\zeta)  e^{-i <z,\zeta>} dz| \leq C_n |\zeta|^{-n}.
\]
\end{lemma}

For small   $|\zeta|\, |\xi|^{-1}$ we  integrate by parts in the variable $\zeta$ and really use
assumptions of Theorem~\ref{AM}.

\begin{lemma}\label{BC2}
For any $\chi\in C_0^\infty ({\Bbb R}^d)$
\[
 |\int_{{\Bbb R}^d}(\partial^\alpha_\xi D^\alpha_{x^\prime}{\bf a})  (x,x+z,\xi
+t\zeta) \chi(\zeta/|\xi|) e^{-i <z,\zeta>} d\zeta| \leq C_n |\xi|^{m+d-|\alpha|(1-r)}.
\]
\end{lemma} 
{\it Proof.} --
We omit the inessential variable $x$ and set $\theta(z,\eta)= \Theta (x,x+z,\eta)$. According to
(\ref{eq:BC2a}), it suffices to check that
\begin{equation}
 \int_{{\Bbb R}^d}\partial^\alpha_\xi \Bigl(e^{i \theta( z,\xi+t\zeta)}  b(z,\xi +t\zeta)\Bigr)
\chi(\zeta/|\xi|) e^{-i <z,\zeta>} d\zeta = O(|\xi|^{p+d-|\alpha|}),
\label{eq:AN2}\end{equation}
if $b\in {\cal S}^p,\: \theta\in {\cal S}^r$ and   $\theta(0,\xi)=0$ for all $\xi$. Estimate
(\ref{eq:AN2}) is, of course, true for $\alpha=0$. We assume it for $|\alpha|=N-1$ and verify for 
$|\alpha|=N$. Note that for some $k$ and $\beta$ with $|\beta|=N-1$
\[
 \partial^\alpha_\xi  (e^{i \theta} b)= \partial^\beta_\xi  (e^{i \theta} b_{\xi_k} +i e^{i \theta}
\theta_{\xi_k} b)
\]
and split up integral (\ref{eq:AN2}) into two terms. Since $b_{\xi_k}\in {\cal S}^{p-1}$, the
integral  containing $\partial^\beta_\xi  (e^{i \theta} b_{\xi_k})$ is bounded, by the inductive
assumption, by $C|\xi|^{p+d-N}$. So we need only to consider 
\begin{equation}
i \int_{{\Bbb R}^d}\partial^\beta_\xi \Bigl(e^{i \theta( z,\xi+t\zeta)}  f (z,\xi +t\zeta)\Bigr)
\chi(\zeta/|\xi|) e^{-i <z,\zeta>} d\zeta ,
\quad f = b  \theta_{\xi_k}\in {\cal S}^{p-1+r},\quad f (0,\xi) =0.
\label{eq:BC15}\end{equation}
Integrating here by parts we find that this integral equals
\[
 \sum_{j=1}^d z_j |z|^{-2}\int_{{\Bbb R}^d}
\partial_{\zeta_j} \{\partial_\xi^\beta \Bigl(   e^{i \theta( z,\xi+t\zeta)} f  ( z,\xi +t\zeta)
\Bigr)\chi (\zeta/|\xi|) \} e^{-i <z,\zeta>}  d\zeta.
\]
Note that 
\begin{equation}
\partial_{\zeta_j} \{\partial_\xi^\beta (   e^{i \theta } f )\chi  \}=|\xi|^{-1} \partial_\xi^\beta
(   e^{i \theta } f )\chi_{\zeta_j}+t \partial_\xi^\beta (   e^{i \theta } f_{\xi_j} )\chi
+ it \partial_\xi^\beta (   e^{i \theta } \theta_{\xi_j}f )\chi.
\label{eq:AN4}\end{equation}
Since
\[
z_j |z|^{-2} f(z,\xi ) \in {\cal S}^{p-1+r} \subset {\cal S}^p,
\]
by our inductive assumption, the integrals containing the first two terms in the right-hand side of 
(\ref{eq:AN4}) are bounded by $C|\xi|^{p+d-N}$. So it remains to consider the integrals
\[
\int_{{\Bbb R}^d}\partial^\beta_\xi \Bigl(e^{i \theta( z,\xi+t\zeta)} g_{j}(z,\xi +t\zeta)\Bigr)
\chi(\zeta/|\xi|) e^{-i <z,\zeta>} d\zeta,
\quad  g_{j}=z_j |z|^{-2} \theta_{\xi_j} f\quad
g_{j}(0,\xi)=0.
\]
They  have the same form as  (\ref{eq:BC15}) but $g_{j}\in  {\cal S}^{p-2+2r}$. Therefore we can
repeat the arguments above. After $n$ steps we arrive at integral  (\ref{eq:BC15}) with
$f \in {\cal S}^q $, where $q=p-\min\{1,  (n+1)(1-r)\}$. So one needs only to choose $n$ such
that  $ (n+1)(1-r)>1.\quad\Box$

Combining Lemmas~\ref{AM11} and \ref{BC2} we obtain estimate (\ref{eq:AM3}) with $q=m+d$. As was
already mentioned, Lemma~\ref{AM1} directly implies Theorem~\ref{AM}.

{\bf 2.3.}
 In this subsection we consider   PDO $A_1,A_2$  defined by formulas
(\ref{eq:BC1}), (\ref{eq:BC2}). Our goal is to
 obtain  for symbols of the operators
$A_1 A_2^\ast  $ and $A_2^\ast A_1 $ the same expansions as in the case $a_j\in {\cal
S}^m_{\rho,\delta}$ with $\rho>\delta$.
 In view of applications to   scattering theory we  suppose that
symbols $a_1$ and $a_2$ are compactly supported in $x$,
 that is for some compact $K_0\subset X$
\begin{equation}
  a_j(x,\xi)=0 \quad{\rm if}\quad x\not\in K_0, \;\forall \xi\in{\Bbb R}^d, \; j=1,2.
\label{eq:CS}\end{equation}
 Then we may assume that $X={\Bbb R}^d$. In this case $A$ and
$A^\ast$ send the Schwartz space $\cal S$ into itself so that the products 
$A_1 A_2^\ast  $ and $A_2^\ast A_1 $ are correctly defined on $\cal S$. Note also that
Theorems~\ref{BC1} and \ref{bc1} hold as well true if PDO $A_1$ and $A_2$ are properly supported.

 The results on the  operator $A_1 A_2^\ast  $ are summarized in the following

\begin{theorem}\label{BC1}
  Suppose that $A_j\in   {\cal C}^{m_j}(\Phi)$ for $j=1,2$ and some
numbers $m_j$. Then $G=A_1 A_2^\ast$ is a PDO with symbol $g\in {\cal S}^m$ for $m=m_1+m_2$ and 
$g(x,\xi)$ admits the asymptotic expansion
\begin{equation}
  g(x,\xi)=\sum_{|\alpha|\geq 0} (\alpha !)^{-1} g_\alpha (x,\xi),
\quad {\rm where}\quad
  g_\alpha(x,\xi)=\partial^\alpha_\xi (a_1(x,\xi) \overline{D^\alpha_x a_2(x,\xi)});
\label{eq:BC6}\end{equation}
 in particular, $g_\alpha\in{\cal S}^{m-|\alpha|(1-r)}$ for all $\alpha$.
\end{theorem}
{\it Proof.} --
The operator $A_2^\ast$ is a PDO with the amplitude $\overline{a_2(x^\prime,\xi)}$ and,
consequently,
\[
(\widehat{A_2^\ast u})(\xi)=(2\pi)^{-d/2}\int_{{\Bbb R}^d}e^{-i<x^\prime,\xi>}
\overline{a_2(x^\prime,\xi)} u(x^\prime) dx^\prime.
\]
Comparing this expression with definition (\ref{eq:BC1}) of the PDO $A_1$ we see that the product
$G=A_1 A_2^\ast$ can be written in the form (\ref{eq:RH1}) with the amplitude
\[
  {\bf g}(x,x^\prime,\xi)=   a_1(x,\xi) \overline{a_2(x^\prime,\xi)}.
\]
Taking into account representations (\ref{eq:BC2}) for symbols $a_j$ we see that
\[
  {\bf g}(x,x^\prime,\xi)=     e^{i\Theta  (x,x^\prime,\xi) 
} b_1(x,\xi) \overline{b_2(x^\prime,\xi)},\quad {\rm where} \quad \Theta 
(x,x^\prime,\xi)=\Phi(x,\xi)- \Phi(x^\prime,\xi).
\]
This amplitude satisfies the assumptions of Theorem~\ref{AM} so that $g \in {\cal S}^m$ and
expansion (\ref{eq:BC6}) is a particular case of (\ref{eq:AM1}).$\quad\Box$

The results on the operator $ A_2^\ast A_1 $ are formulated similarly to Theorem~\ref{BC1}
although their proof is somewhat different.

\begin{theorem}\label{bc1}
  Suppose that $A_j\in   {\cal C}^{m_j}(\Phi)$ for $j=1,2$ and some
numbers $m_j$. Then $H=A_2^\ast A_1 $ is a PDO with symbol  $h\in{\cal S}^m $ for $m=m_1+m_2$ and
$h(x,\xi)$ admits the asymptotic expansion
\begin{equation}
  h(x,\xi)=\sum_{|\alpha|\geq 0} (\alpha !)^{-1} h_\alpha (x,\xi),
\quad {\rm where}\quad
  h_\alpha(x,\xi)=D^\alpha_x (a_1(x,\xi) \overline{\partial^\alpha_\xi a_2(x,\xi)});
\label{eq:bc2}\end{equation}
in particular, $h_\alpha\in{\cal S}^{m-|\alpha|(1-r)}$ for all $\alpha$.
\end{theorem}

Note that $h_\alpha\in{\cal S}^{m-|\alpha|(1-r)}$ because the oscillating factors $e^{\pm i \Phi 
(x,\xi) }$ cancel each other in the product $a_1 \overline{\partial^\alpha_\xi a_2}$. 

Using again that $A_2^\ast$ is the PDO with   amplitude $\overline{a_2(x^\prime,\xi)}$ 
we obtain the equality
\[ 
(A_2^\ast A_1 u)(x)=(2\pi)^{-3d/2}\int_{{\Bbb R}^d}\int_{{\Bbb R}^d}e^{i<x-x^\prime,\eta>}
\overline{a_2(x^\prime,\eta)} \Bigl(
\int_{{\Bbb R}^d}
a_1(x^\prime,\xi) e^{i<x^\prime,\xi>} \hat{u}(\xi) d \xi\Bigr) dx^\prime d\eta.
\]
Thus, $H=A_2^\ast A_1 $   is the PDO with  symbol
\begin{equation}
  h(x,\xi)= (2\pi)^{-d} \int_{{\Bbb R}^d} \int_{{\Bbb R}^d} a_1(y,\xi)
\overline{a_2(y,\xi+\zeta)} e^{i <x-y,\zeta> } dy d\zeta.
\label{eq:BC7A}\end{equation}
Let us first obtain a rough estimate on this function.

\begin{lemma}\label{ha1}
  Suppose that $a_j\in   {\cal S}^{m_j}_{\rho,\delta}$, $j=1,2$ 
 for some $\rho\geq 0$, $\delta <1$.
Then for all multi-indices $\alpha,\beta$  the function
$\partial_\xi^\alpha \partial_x^\beta h(x,\xi)$
is bounded by  $ C (1+|\xi|)^n $ for some  $n=n(\alpha,\beta)$ and $C=C( \alpha,\beta)$.
\end{lemma}
{\it Proof.} --
 Using (\ref{eq:IP})  and integrating by parts in (\ref{eq:BC7A}) in the variable $y$ we see that
\[
|\partial_\xi^\alpha \partial_x^\beta h(x,\xi)|
\leq 
C  \langle \xi \rangle ^{m_1} \int_{{\Bbb R}^d} \langle \xi+\zeta\rangle ^{m_2} \langle
\zeta\rangle^{|\beta|-k} (\langle
\xi+\zeta\rangle ^{\delta k} +
\langle \zeta\rangle ^{\delta k}) d\zeta.
\]
This gives the necessary estimate with, for example, $n=m_1+m_2^+ +\delta k$ (here
$m_2^+=\max\{0,m_2\}$) 
 if $k$ is large enough, say, $(1-\delta)k >|\beta| +d +m_2^+$. $\quad\Box$

For any fixed $x,\xi$, expression (\ref{eq:BC7A}) coincides with the left-hand side of
(\ref{eq:Ta2})   where
\[
 {\bf p}(z,\zeta;x,\xi) =  a_1( x+z,\xi ) \overline{ a_2( x+z,\xi +\zeta )}.
\]
 Then the coefficients $(\partial^\alpha_\zeta D^\alpha_z  {\bf p})  (0,0;x,\xi)$ equal to the
numbers $h_\alpha(x,\xi)$ defined in (\ref{eq:bc2}). Equality (\ref{eq:Ta2}) shows that $h(x,\xi)$ is
a sum of the terms $(\alpha!)^{-1} h_\alpha(x,\xi)$ over $\alpha,\; 0\leq |\alpha|\leq N-1$, and of
the rest
$h^{(N)}(x,\xi)$ defined by (\ref{eq:Ta3}). Thus, for the proof of Theorem~\ref{bc1}  we need only to
estimate $h^{(N)}(x,\xi)$. More precisely, it suffices (cf. the proof of Theorem~\ref{AM}) to obtain
estimates (\ref{eq:AM3}), uniformly in $t\in [0,1]$, for the functions
\begin{equation}
  R_\alpha(x,\xi;t)= \int_{{\Bbb R}^d}  v_\alpha ( \xi,\zeta;t) e^{i <x,\zeta> }   d\zeta,
\label{eq:bc1}\end{equation}
where
\begin{equation}
 v_\alpha (\xi,\zeta;t)=  \zeta^\alpha \int_{{\Bbb R}^d}  a_1(y,\xi)
\overline{(\partial^\alpha_\xi a_2)(y,\xi+t\zeta)} e^{-i <y,\zeta> }   dy.
\label{eq:bc1a}\end{equation}

 Remark first that, quite similarly to Lemma~\ref{ha1}, one can check that under  its assumptions,
for any $n$,
\begin{equation}
 | v_\alpha ( \xi,\zeta;t)| \leq C_n (1+|\zeta|)^{-n}
\quad {\rm if}\quad |\zeta|\geq |\xi|/2.
\label{eq:bc1b}\end{equation}
An estimate of function (\ref{eq:bc1a}) for $|\zeta|\leq |\xi|/2$ requires condition
(\ref{eq:BC2}).

\begin{lemma}\label{ha2}
Under the assumptions of Theorem~$\ref{bc1}$,
\begin{equation}
 | v_\alpha (\xi,\zeta;t)| \leq C  (1+|\xi|)^{m-|\alpha|(1-r)}
\quad {\rm if}\quad |\zeta|\leq |\xi|/2.
\label{eq:bc1bx}\end{equation}
\end{lemma}

Estimating integral (\ref{eq:bc1}) with the help of (\ref{eq:bc1b}) and (\ref{eq:bc1bx}) we obtain
bound (\ref{eq:AM3}), where $q=m+d$, for  $R_\alpha(x,\xi;t)$.  So to  conclude the proof of
Theorem~\ref{bc1} it remains to verify Lemma~\ref{ha2}.

The proof of estimate (\ref{eq:bc1bx}) relies on the representation
\begin{equation}
 a_1(y,\xi) \overline{(\partial^\alpha_\xi
a_2)(y,\eta)}=  e^{i\varphi(y,\xi,\eta)} w_\alpha(y,\xi, \eta),
\label{eq:bc5}\end{equation}
where
\begin{equation}
\varphi (y,\xi, \eta)= \Phi (y,\xi)  - \Phi (y, \eta)
\label{eq:bc6}\end{equation}
and, according to (\ref{eq:BC2}), for all $\beta$  
\begin{equation}
| (\partial_y^\beta w_\alpha)(y,\xi,\xi+t\zeta)| \leq C_\beta 
(1+|\xi|)^{m-|\alpha|(1-r)}.
\label{eq:bc1b1}\end{equation}
 This estimate for $\beta=0$ implies already  that
(\ref{eq:bc1bx}) is fulfilled for $|\zeta|$ bounded. To obtain (\ref{eq:bc1bx}) in the whole ball
$|\zeta|\leq |\xi|/2$ we need to get rid of the growing factor $\zeta^\alpha$ in (\ref{eq:bc1a}).
Note that function (\ref{eq:bc6})  satisfies,  for all $\beta$ and any compact $K $, the inequality 
\begin{equation}
\sup_{y\in K}   |(\partial^\beta_y \varphi) (y,\xi,\xi+t\zeta)|\leq C |\xi|^{r-1}|\zeta|,
\label{eq:bc8}\end{equation}
which is a consequence of the assumption
 $\Phi\in {\cal S}^r$. Let us use the following

\begin{lemma}\label{BC2A}
Let a function $\varphi$ satisfy $(\ref{eq:bc8})$ and let $1\leq |\zeta|\leq |\xi|/2$.
 Suppose that a function $w(y,\xi,\zeta,t)$ is compactly supported in the variable $y$ and  for some
$p$ and all $\beta$
\begin{equation}
 | (\partial_y^\beta w )(y,\xi,\zeta,t)| \leq C  |\xi|^p,\quad
|\zeta|\leq |\xi|/2.
\label{eq:bc1c}\end{equation}
Then
\begin{equation}
  |\zeta| \int_{{\Bbb R}^d} w (y,\xi,\zeta,t) e^{i\varphi(y,\xi,\zeta,t)} e^{-i <y,\zeta> }  
dy = \int_{{\Bbb R}^d} \tilde{ w} (y,\xi,\zeta,t) e^{i\varphi(y,\xi,\zeta,t)} e^{-i
<y,\zeta> }   dy,
\label{eq:bc1aa}\end{equation}
where $\tilde{w}$ is compactly supported in $y$ and  satisfies $(\ref{eq:bc1c})$ for all $\beta$.
\end{lemma}
{\it Proof.} --
Integrating in the left-hand side of (\ref{eq:bc1aa}) by parts we rewrite it as
\[
-i \sum_{j=1}^d \zeta_j |\zeta|^{-1}\int_{{\Bbb R}^d}
\partial_{y_j} \Bigl( w (y,\xi,\zeta,t) e^{i \varphi(y,\xi, \zeta,t)}\Bigr)
  e^{-i <y,\zeta>}  dy.
\]
The functions
$  (\partial_{y_j}  w ) (y,\xi,\zeta,t)$ satisfy   (\ref{eq:bc1c}) and, by
virtue of (\ref{eq:bc8}), the functions
\begin{equation}
 |\zeta|^{-1}  w  (y,\xi,\zeta,t) (\partial_{y_j}  \varphi)
(y,\xi,\zeta,t)
\label{eq:byx}\end{equation}
 satisfy (\ref{eq:bc1c}) with $p$ replaced by $p+r-1$. Let us again integrate  by parts  in
integrals containing functions (\ref{eq:byx}).
  Then after $n   $ steps, we arrive at representation  (\ref{eq:bc1aa}) with a function $\tilde{ w}$
bounded (with its derivatives in $y$) by
\[
C|\xi|^p (1+ |\zeta| |\xi|^{-n(1-r)}).
\]
If $n >(1-r)^{-1} $, this gives inequality (\ref{eq:bc1c}) for  $\tilde{ w}$.
$\quad\Box$

Let us now take into account representation (\ref{eq:bc5}) and apply Lemma~\ref{BC2A}
 $N=|\alpha|$ times to integral (\ref{eq:bc1a}). Thus, the growing factor $\zeta^\alpha$ may be
``eaten up" which, by virtue of (\ref{eq:bc1b1}), gives (\ref{eq:bc1bx}).
 This  conclude the proof of Lemma~\ref{ha2} and hence of
Theorem~\ref{bc1}.

{\bf 2.4.} 
A simple tradional condition of boundedness of PDO is formulated in the following elementary

\begin{lemma}\label{Bc1} 
Suppose that $b\in {\cal S}^m_{0,0}$ and $b(x,\xi)$ is compactly supported
in the variable $x$. Then the PDO with symbol $b(x,\xi)$ is bounded  in the space $L_2({\Bbb R}^d)$
if
$m=0$ and it is compact if $m<0$.
 \end{lemma}

A proof can be easily deduced from the expansion of $b(x,\xi)$ in the Fourier series in the
variable $x$.

Combining this assertion with Theorem~\ref{BC1},  we obtain conditions
of boundedness of PDO with oscillating symbols.  

\begin{theorem}\label{Bc2}
 Let $a\in {\cal C}^m(\Phi)$ and let $a(x,\xi)$ be compactly supported in the variable $x$. 
Then the PDO with symbol $a(x,\xi)$ is bounded in the space $L_2({\Bbb R}^d)$ if $m=0$ and it is
compact if $m<0$.
\end{theorem}
 {\it Proof.} --
 By Theorem~\ref{BC1},  symbol of the PDO
$G=A A^\ast$ belongs to the class $ {\cal S}^{2m}_{1,0}$. So according to
 Lemma~\ref{Bc1}  it is bounded  (compact)  which implies   boundedness (compactness) of the
operators $A^\ast$ and $A$.$\quad\Box$

Applications to scattering theory require to single out the singular part of the operators 
$A_1 A_2^\ast $ and $A_2^\ast A_1 $.

\begin{theorem}\label{Bc3}
 Let  $A_j,\; j=1,2,$ be the PDO  with symbols  $a_j\in {\cal C}^0(\Phi)$
$($so that 
$a_j=e^{i\Phi}b_j$ with $b_j\in {\cal S}^0$$)$ satisfying condition $(\ref{eq:CS})$. Denote by
$B$ the PDO with symbol
\[ b(x,\xi)=b_1(x,\xi)\: \overline{b_2(x,\xi)},\quad b\in {\cal S}^0.
\] Then both operators $A_1 A_2^\ast -B $ and $A_2^\ast A_1 -B $  are compact in the space
$L_2({\Bbb R}^d)$.
\end{theorem}
 {\it Proof.} --
 It follows from Theorem~\ref{BC1} that $G_1=A_1 A_2^\ast -B $ is the
PDO with symbol $g_1$ from the class ${\cal S}^{-1+r}$. Since, moreover, $g_1$ is compactly
supported in $x$, the operator $G_1$ is compact by virtue of Lemma~\ref{Bc1}.

Similarly, by Theorem~\ref{bc1}, $H_1=A_2^\ast A_1 -B $ is the PDO with  symbol $h_1$ from the class
${\cal S}^{-1+r}$. Let $\Omega$ be multiplication by a function $\Omega\in C_0^	\infty ({\Bbb
R}^d)$. By virtue of Lemma~\ref{Bc1}, the operator $\Omega H_1$ is compact. Choose
$\Omega=\bar{\Omega}$ such that $\Omega(x)=1$ in a neighbourhood of $K_0$. Then $(I-\Omega) B=0$.
Let us check that the operator $(I-\Omega) A^\ast_2$ is compact. 
Clearly, $A_2$ is an integral operator with kernel 
$k_2(x, y)={\bf k}_2(x,x-y),$ where ${\bf k}_2(x,z)$ is a $C^\infty$-function outside the diagonal
$z=0$; moreover,
${\bf k}_2(x,z)=0$ if $x\not\in K_0$ and ${\bf k}_2(x,z)\rightarrow 0$
 quicker than any power of $|z|^{-1}$ as $|z|\rightarrow\infty$. Therefore the kernel
$ {\bf k}_2(x,x-y) (1-\Omega(y))$ of the operator $ A_2 (I-\Omega)$ is a
$C^\infty$-function which is compactly supported in $x$ and is rapidly decreasing as
$|y|\rightarrow\infty$. This implies compactness of $(I-\Omega)  A_2^\ast$ and hence of
$(I-\Omega)  H_1$.$\quad\Box$

{\bf 2.5.} 
Here we use the stationary phase method to show that a PDO $A$  defined by   oscillating 
amplitude ${\bf a}$ admits representation (\ref{eq:BC1})  with   oscillating symbol $a$. Thus we
suppose that $A$ is given by formula (\ref{eq:RH1}),  where   amplitude ${\bf a}$ satisfies
(\ref{eq:I1}). Assume that  a PDO $A$ is properly supported.
 We shall check that   symbol $a$ of the PDO $A$ satisfies (\ref{eq:BC2}) and  find   expressions
for the functions $\Phi$ and $b$ in terms of the functions $\Theta$ and ${\bf b}$. 

Let us proceed from  representation (\ref{eq:BC7}). Inserting there (\ref{eq:I1}),
making   the change of variables $\zeta=|\xi|\eta$ and denoting
\begin{equation}
\Xi(z,\eta; x,\xi)= |\xi|^{-1}\Theta(x,x+z,\xi+ |\xi|\eta) - <z,\eta> 
\label{eq:SP1a}\end{equation}
  we see that
\begin{equation}
  a(x,\xi)= (2\pi)^{-d}|\xi|^d \int_{{\Bbb R}^d} \int_{{\Bbb R}^d}
\exp\Bigl(i|\xi|\Xi(z,\eta; x,\xi)\Bigr) {\bf b}(x,x+z,\xi+|\xi|\eta)  dzd\eta.
\label{eq:SP3}\end{equation}
 The asymptotics of this integral
  as $|\xi|\rightarrow\infty$ is determined by stationary points
$z_s=z_s(x,\xi),\:\eta_s=\eta_s(x,\xi)$ satisfying the equations
\[
\Xi_\eta(z_s,\eta_s; x,\xi)=0,\quad \Xi_z(z_s,\eta_s; x,\xi)=0
\] 
or, in view of (\ref {eq:SP1a}),
\begin{equation} 
z_s=\Theta_\xi (x,x+z_s,\xi+ |\xi|\eta_s),\quad \eta_s= |\xi|^{-1} \Theta_y
(x,x+z_s,\xi+|\xi|\eta_s),
\label{eq:SP2}\end{equation}
 where   $\Theta_y $ and  $\Theta_\xi $ are
derivatives of the function
 $\Theta (x,y,\xi)$ in the second and third variables.

\begin{lemma}\label{SP1}
 For any $x$ and sufficiently large $|\xi|$ system
$(\ref{eq:SP2})$ has a solution $z_s=z_s(x,\xi),\eta_s=\eta_s(x,\xi)$, which can be obtained by
iterations
\[
 z_s=\lim_{p\rightarrow\infty}z_s^{(p)},\quad
\eta_s=\lim_{p\rightarrow\infty}\eta_s^{(p)},
\]
 where $ z_s^{(0)}=0,\, \eta_s^{(0)}=0$,
\[
 z_s^{(p)}=\Theta_\xi(x,x+z_s^{(p-1)},\xi+ |\xi| \eta_s^{(p-1)}), \quad
 \eta_s^{(p)}= |\xi|^{-1} \Theta_y(x,x+z_s^{(p-1)},\xi+ |\xi| \eta_s^{(p-1)}).
\]
 The vector-functions $z_s(x,\xi)$ and $\eta_s(x,\xi)$ are
infinitely differentiable in the variables $x$ and $\xi$ and belong to the class ${\cal S}^{-1+r}$.
\end{lemma}

We omit an elementary proof which relies on the estimate 
\[
| z_s^{(p+1)}-z_s^{(p)}| + | \eta_s^{(p+1)}-\eta_s^{(p)}|
\leq
 C|\xi|^{-1+r} (|z_s^{(p)}-z_s^{(p-1)}|+ |\eta_s^{(p)}-\eta_s^{(p-1)}|).
\]
In its turn, this estimate is a consequence of the assumption $\Theta\in {\cal S}^r$.

Let us introduce the Hessian of the function $\Xi(z,\eta;x,\xi)$ in the variables $z$ and
$\eta$:
\begin{equation}
 h(z,\eta;x,\xi)={\rm Hess}_{z,\eta}\,\Xi(z,\eta;x,\xi) =\left( \begin{array}{cc}
  \partial_{zz}\Xi &  \partial_{z\eta}\Xi\\ 
\partial_{\eta z}\Xi &  \partial_{\eta\eta}\Xi
 \end{array}\right)
\label{eq:SP9}\end{equation}
 and let $\det h$ and ${\rm sgn}\, h$ be the determinant and the signature of this
$(2d)\times(2d)$-matrix. According to (\ref {eq:SP1a}),
\[
\Xi_{zz}(z,\eta;x,\xi )= |\xi|^{-1}\Theta_{yy}(x,x+z,\xi+|\xi|\eta)= O(|\xi|^{-1+r}),
\]
\[
\Xi_{z,\eta}(z,\eta;x,\xi )=-I+ \Theta_{y\xi}(x,x+z,\xi+|\xi|\eta)= -I+ O(|\xi|^{-1+r}),
\]
\[
\Xi_{\eta\eta}(z,\eta;x,\xi)= |\xi| \Theta_{\xi\xi}(x,x+z,\xi+|\xi|\eta)= O(|\xi|^{-1+r}),
\] 
 and, consequently,
\[
 |\det h( z,\eta;x,\xi)|=1+ O(|\xi|^{-1+r}), \quad {\rm sgn}\; h(z,\eta;x,\xi)=0
\]
 for sufficiently large $|\xi|$.

Applying the stationary phase method to integral (\ref{eq:SP3}) we obtain

\begin{theorem}\label{SP4}
 Suppose that a PDO $A$ is given by formula $(\ref{eq:RH1})$,  where the
amplitude ${\bf a}$ admits  representation $(\ref{eq:I1})$. Then   symbol of the operator $A$
satisfies
$(\ref{eq:BC2})$, where
\[
\Phi(x,\xi)=\Theta(x,x+z_s,\xi+ |\xi| \eta_s) - |\xi| <z_s,\eta_s>;
\]
 in particular, $\Phi\in{\cal S}^r$ and
\begin{equation}
 \Phi(x,\xi)=\Theta(x,x,\xi) + \Phi_0(x,\xi),\quad {\rm where }\quad \Phi_0 \in{\cal S}^{-1+2r}.
\label{eq:SP18}\end{equation}
 The function $b\in{\cal S}^m$ and
\[
 b(x,\xi)={\bf b}(x,x,\xi) +O(|\xi|^{m-1+r}),\quad |\xi|\rightarrow\infty.
\]
\end{theorem}

This result is, of course, of interest only in the case $r\geq 1/2$ when the standard expression
(see e. g. \cite{Sh}) of the symbol in terms of the amplitude is not applicable and the contribution
of $\Phi_0(x,\xi)$ in (\ref{eq:SP18}) cannot be neglected.
 
Combining  Theorems~\ref{Bc2} and \ref{SP4},  we obtain conditions of boundedness of PDO with
oscillating amplitudes.

\begin{theorem}\label{Bc2d}
 Let  ${\bf a}\in {\cal C}^m(\Theta)$ and let ${\bf a}(x,x^\prime,\xi)$ be
compactly supported in the variables $x$ and $x^\prime$.  Then the PDO with   amplitude ${\bf
a}(x,x^\prime,\xi)$ is bounded
$($compact$)$ in the space $L_2({\Bbb R}^d)$ if $m=0$ $(m<0)$.
\end{theorem}
 
\section {the action on an exponent}

{\bf 3.1.}
Our  goal here is to calculate a PDO $A$ defined by (\ref{eq:BC1}), (\ref{eq:BC2}) on functions
\begin{equation}
 u_{\lambda,\varepsilon}(x)=e^{i\lambda\psi(x )} f_\varepsilon(x),\quad {\rm where}\quad
f_\varepsilon(x)=\varepsilon^{-d/2}f((x-x_0)/\varepsilon),
\label{eq:I2d}\end{equation}
 $f\in C_0^\infty (X)$, $\varepsilon\in(0,\varepsilon_0)$  and $\lambda\rightarrow\infty$.
As in subsection 2.5, we use the stationary phase  method.
We allow $\varepsilon\rightarrow 0$ when $f=f_\varepsilon$ are shrinking to the  point $x_0$.
All our estimates will be uniform  with respect to phase functions $\psi\in C^\infty$ satisfying
the following

\begin{assumption}\label{psi}
For a neighbourhood $U$ of the point $x_0$ and constants $C_\alpha$ and $c$,
\[ 
|\partial^\alpha \psi(x)|\leq C_\alpha  \quad {\rm and} \quad
|\psi^\prime(x )|\geq c>0, \quad x\in U.
\] 
\end{assumption}

If $\varepsilon\rightarrow0$, then automatically ${\rm supp}\,f_\varepsilon \subset U$ for
sufficiently small
$\varepsilon$. If $\varepsilon$ is fixed, then we require that ${\rm supp}\,f_\varepsilon \subset
U$.

Let us insert (\ref{eq:I2d}) into (\ref{eq:RH1}), where ${\bf a}(x,x^\prime,\xi)=a(x,\xi)$.  Making
the change of variables $\xi=\lambda\eta$, we see that
\begin{equation}
 (Au_{\lambda,\varepsilon})(x)=(2\pi)^{-d}\lambda^d e^{i\lambda \psi(x)}
\int_{X} \int_{{\Bbb R}^d}e^{i\Gamma(y,\eta;x,\lambda)} b(x,\lambda\eta) f_\varepsilon(y) dy d\eta,
\label{eq:S2}\end{equation}
where
\begin{equation}
\Gamma(y,\eta;x,\lambda)
=\lambda < x-y, \eta > +  \Phi(x,\lambda \eta)+ \lambda\psi (y )-                        \lambda
\psi(x).
\label{eq:S3}\end{equation}
The  critical points $y_s=y_s(x,\lambda), \eta_s=\eta_s(x,\lambda)$ of this function are
defined by the equations
\[
 \Gamma_\eta(y_s,\eta_s;x,\lambda) =0,\quad \Gamma_y(y_s,\eta_s;x,\lambda) =0
\]
or, in view of (\ref{eq:S3}),
\begin{equation}
y_s=x+\Phi_\xi(x,\lambda \eta_s),\quad \eta_s=\psi^\prime (y_s ).
\label{eq:S5}\end{equation}
This gives an equation for $y_s$:
\begin{equation}
 y_s=x+ \Phi_\xi(x ,\lambda \psi^\prime (y_s)).
\label{eq:S5a}\end{equation}
We omit the proof of the following elementary assertion.

\begin{lemma}\label{S1}
 Let $x\in X$ and $\psi^\prime(x)\neq 0$. Then for sufficiently large
$\lambda$ equation
$(\ref{eq:S5a})$ has a unique solution $y_s=y_s(x,\lambda)$, which can be obtained by iterations
\[
 y_s=\lim_{p\rightarrow\infty}y_s^{(p)},
\]
 where $ y_s^{(0)}=x$,
\[
  y_s^{(p)}=x+\Phi_\xi(x,\lambda \psi^\prime (y_s^{(p-1)} )) .
\]
 In particular,
\begin{equation}
 y_s(x,\lambda)=x+O(\lambda^{-1+r}).
\label{eq:S8}\end{equation}
Furthermore, the vector-function $y_s(x,\lambda)$ is infinitely differentiable in the variable $x$
and
\begin{equation}
 y_s^\prime(x,\lambda)=I+O(\lambda^{-1+r}),\quad 
y_s^{(\alpha)}(x,\lambda)=O(\lambda^{-1+r}),\; |\alpha| \geq 2.
\label{eq:S8df}\end{equation}
\end{lemma}

It follows now from the second equation (\ref{eq:S5}) and (\ref{eq:S8}) that
\begin{equation}
\eta_s(x,\lambda)=\psi^\prime(x  )+O(\lambda^{-1+r}).
\label{eq:S8a}\end{equation}
Let us also introduce the phase $\Gamma$ at the critical point:
\begin{equation}
G(x,\lambda)=\Gamma(y_s(x,\lambda), \eta_s(x,\lambda);x,\lambda).
\label{eq:GG}\end{equation}

It is convenient to transfer  in (\ref{eq:S2}) the dependence of $f_\varepsilon$ on $\varepsilon$
into the phase function.
Setting $x=x_0+\varepsilon w$ and making the change of variables $y=x_0+\varepsilon z$ we obtain that
\begin{equation}
 (Au_{\lambda,\varepsilon})(x_0+\varepsilon w)=(2\pi)^{-d}\lambda^d \varepsilon^{d/2}
\int_{X} \int_{{\Bbb R}^d}e^{i\lambda\varepsilon\Xi(z,\eta;w,\lambda,\varepsilon)} b(x_0+\varepsilon
w,\lambda\eta)
 f(z) dz d\eta,
\label{eq:S2A}\end{equation}
 where
\begin{equation}
\Xi(z,\eta;w,\lambda,\varepsilon)=(\lambda\varepsilon)^{-1}\Gamma(x_0+\varepsilon
z,\eta;x_0+\varepsilon w,\lambda).
\label{eq:S3A}\end{equation}
The stationary phase method can be applied to integral (\ref{eq:S2A}) if
 $\varepsilon \lambda^{1-r}\rightarrow\infty$. Indeed, 
its stationary point $z_s,\eta_s$ are given, according to
(\ref{eq:S8}), (\ref{eq:S8a}), by relations
\begin{equation}
z_s=\varepsilon^{-1}(y_s(x_0+\varepsilon w,\lambda)-x_0)=w+O(\varepsilon^{-1}\lambda^{-1+r}),
\label{eq:S5b}\end{equation}
\begin{equation} 
\eta_s=\eta_s(x_0+\varepsilon w,\lambda)=\psi^\prime(x_0+\varepsilon w  )+O(\lambda^{-1+r})
\label{eq:S5c}\end{equation}
and consequently have finite limits as $\varepsilon^{-1}\lambda^{-1+r}\rightarrow 0$. Furthermore,
it follows from (\ref{eq:S3}), (\ref{eq:S3A}) that
\begin{equation}
\Xi_{zz} =\varepsilon\psi^{\prime\prime}(x_0+\varepsilon z),
\quad  \Xi_{z\eta} =-I,\quad
\Xi_{\eta\eta} =\lambda \varepsilon^{-1} 
\Phi_{\xi\xi}(x_0+\varepsilon w,\lambda\eta)=  O(\varepsilon^{-1}\lambda^{-1+r}),
\label{eq:S12}\end{equation}
\[
 (\partial^\alpha_{z}\Xi)(z,\eta;w,\lambda,\varepsilon)=O(\varepsilon^{|\alpha|-1}),
\quad (\partial^\alpha_{\eta}\Xi)(z,\eta;w,\lambda,\varepsilon)= 
O(\varepsilon^{-1}\lambda^{-1+r}),\quad |\alpha|\geq 2,
\]
 and mixed derivatives of $\Xi$ of order higher than two are zero. Let
\[
 h (z,\eta;w,\lambda,\varepsilon)={\rm Hess}_{z,\eta}\,\Xi(z,\eta;w,\lambda,\varepsilon)
\]
be the Hessian (cf. (\ref{eq:SP9})) of the function $
\Xi(z,\eta;w,\lambda,\varepsilon)$ in the variables $z$ and $\eta$.
By virtue of (\ref{eq:S12}),  the determinant  and the signature
 of this $(2d)\times(2d)$-matrix satisfy the relations
\[
 |\det h(z,\eta;w,\lambda,\varepsilon)|=1+ O(\lambda^{-1+r}), \quad  {\rm sgn}\;
h(z,\eta;w,\lambda,\varepsilon)=0
\]
 as $\lambda\rightarrow\infty$ uniformly in  $\varepsilon\in
(0,\varepsilon_0)$. Thus, the stationary phase method gives
the following   result for integral (\ref{eq:S2A}).
 
\begin{lemma}\label{S2} 
If $\varepsilon\lambda^{1-r}\rightarrow\infty$, then
\begin{eqnarray}
\varepsilon^{d/2} (Au_{\lambda,\varepsilon})(x_0+\varepsilon w)
=   e^{i\lambda\varepsilon \Xi(z_s,\eta_s;w,\lambda,\varepsilon)} e^{i\lambda \psi(x)}
\nonumber\\
\times b(x_0+\varepsilon w,\lambda\eta_s) f(z_s) \Bigl(1+O(\lambda^{-1+r})\Bigr) +
O(\varepsilon^{-1}\lambda^{m-1})
\label{eq:S10}\end{eqnarray}
uniformly in $w,\; |w| \leq c$. Moreover,   if $ f(x)=0$ for   $ |x| \geq c_0$, then for any $n$
\begin{equation}
\varepsilon^{d/2} (Au_{\lambda,\varepsilon})(x_0+\varepsilon w) =  O(\lambda^m
(\lambda \varepsilon)^{-n}),\quad |w| \geq c_1> c_0.
\label{eq:S11}\end{equation}
\end{lemma}

Note that the second assertion can be proven by integration by parts in (\ref{eq:S2A}) in the
variable $\eta$ since
\[
\Xi_\eta (z,\eta;w,\lambda,\varepsilon) =z-w+\varepsilon^{-1}\Phi_\eta (x_0+\varepsilon
z,\lambda\eta) =z-w+ O(\varepsilon^{-1} \lambda ^{-1+r})
\]
 is separated from zero.

The right-hand side of (\ref{eq:S10}) can be simplified. By virtue of (\ref{eq:GG}), (\ref{eq:S3A}),
\[
\lambda \varepsilon\Xi(z_s,\eta_s; w,\lambda,\varepsilon)= G(x,\lambda),
\quad x=x_0+\varepsilon w,
\]
and, by virtue of (\ref{eq:S5b}), (\ref{eq:S5c}),
\begin{equation}
  f(z_s)=f(w)+O(\varepsilon^{-1} \lambda ^{-1+r}),
\label{eq:S13a}\end{equation}
\begin{equation}
 b(x ,\lambda\eta_s)= b(x,\lambda\psi^\prime (x )) +O(\lambda^{m-1+r}).
\label{eq:S13}\end{equation}
Now we can check relation (\ref{eq:I3}). 

\begin{proposition}\label{S5}
Suppose that   symbol $a(x,\xi)$ of  PDO $(\ref{eq:BC1})$ satisfies conditions $(\ref{eq:BC2})$.
 Let functions $u_{\lambda,\varepsilon}$ be defined by $(\ref{eq:I2d})$,  let
$K\subset X$ be any compact set and let the function
$G(x,\lambda)$ be defined by equations $(\ref{eq:S5})$ and equalities $(\ref{eq:S3})$,
$(\ref{eq:GG})$.  Then 
\begin{equation}
 (Au_{\lambda,\varepsilon})(x ) =  e^{iG(x,\lambda)} 
 b(x ,\lambda \psi^\prime(x )) u_{\lambda,\varepsilon}(x)  +
R_{\lambda,\varepsilon}(x ),
\label{eq:I3A}\end{equation}
where the norm of $R_{\lambda,\varepsilon}$  in
$L_2(K)$ is bounded by $C\varepsilon^{-1} \lambda^{m-1+r}$ as   $ \varepsilon\lambda^{1-r}
\rightarrow\infty$ uniformly    with respect to functions $\psi$ satisfying Assumption~$\ref{psi}$.
\end{proposition}
{\it Proof.} --
Let the function  $R_{\lambda,\varepsilon}$ be defined by equality  (\ref{eq:I3A}) and let $c$ be
some positive number. Suppose that $ f(x)=0$ for  $ |x| \geq c_0$ and choose $c_1>c_0$.
Making the change of variables $x=x_0+ \varepsilon w$ we see that it suffices to
estimate 
\[
 \int_{|x-x_0|\leq  c} |R_{\lambda,\varepsilon}(x)|^2 dx= \varepsilon^d
\int_{|w|\leq  c_1} |R_{\lambda,\varepsilon}(x_0+\varepsilon w)|^2 dw
+\varepsilon^d
\int_{c_1 \leq |w|\leq  c \varepsilon^{-1}} |R_{\lambda,\varepsilon}(x_0+\varepsilon w)|^2 dw.
\]
According to (\ref{eq:S10}),  (\ref{eq:S13a}) and (\ref{eq:S13}), the function
$\varepsilon^{d/2} R_{\lambda,\varepsilon}(x_0+\varepsilon w)$ is bounded by $C\varepsilon^{-1}
\lambda^{m-1+r}$, which gives the necessary estimate for the first integral in the right-hand side.
It follows from (\ref{eq:S11}) that the second
integral  is bounded by $C_1 \varepsilon^{-d}\lambda^{2m}(\varepsilon
\lambda)^{-2n}$. Since  $n$ is arbitrary and  $\varepsilon\lambda^{1-r}\rightarrow\infty$, this term
can be estimated by  $C\varepsilon^{-2} \lambda^{2(m-1+r)}$.
$\quad\Box$

\noindent{\bf Remark}$\;$
Since
\begin{equation}
b(x,\lambda\psi^\prime (x ))=b(x_0,\lambda\psi^\prime (x_0 )) +O(\lambda^m|x-x_0|),
\label{eq:I3An}\end{equation}
the function $b(x,\lambda\psi^\prime (x ))$ in (\ref{eq:I3A}) can be
replaced by the number $b(x_0,\lambda\psi^\prime (x_0 ))$ if $\varepsilon\rightarrow 0$.
In this case  
$||R_{\lambda,\varepsilon}||_{L_2(K)}$ is estimated by $C\lambda^m(\varepsilon^{-1}
\lambda^{-1+r}+ \varepsilon)$.

\noindent{\bf Remark}$\;$
 Of course, Proposition~\ref{S5} holds for phase functions $\psi(x,\lambda)$ depending on the
parameter $\lambda$ as long as Assumption~\ref{psi} is satisfied.

{\bf 3.2.}
In the next section we shall need  the asymptotics for large $\lambda$ of the function $
G(x,\lambda)$ constructed in Proposition~\ref{S5}. Let us show that its leading term  is $\Phi
(x,\lambda \psi^\prime (x))$ which grows as $\lambda^r$.  Set
\begin{equation}
\Omega_1 (x,\lambda)=\lambda^{1-2r}\Bigl(\Phi (x,\lambda \psi^\prime (y_s ))-\Phi (x, 
\lambda\psi^\prime (x ))\Bigr),
\label{eq:Om1}\end{equation}
\begin{equation}
\Omega_2 (x,\lambda)= \lambda^{2-2r} \Bigl(\psi( y_s )- 
\psi(x ) - < y_s-x,\psi^\prime(y_s )>\Bigr)
\label{eq:Om2}\end{equation}
and
\begin{equation}
\Omega (x,\lambda)= \Omega_1 (x,\lambda)+\Omega_2 (x,\lambda).
\label{eq:Om}\end{equation}
Then it follows from (\ref{eq:S3})  that
\begin{equation}
 G(x,\lambda)=\Phi (x,\lambda \psi^\prime (x ))+ \lambda^{-1+2r} \Omega (x,\lambda).
\label{eq:S16}\end{equation}
 To estimate $\Omega$ it is convenient to rewrite  (\ref{eq:Om1}), (\ref{eq:Om2}) in the following
form.

\begin{lemma}\label{om} 
Set $ Y_s(x,t)= ty_s+(1-t)x$. Then
\begin{equation}
\Omega_1 (x,\lambda)=\lambda^{2-2r}\int_0^1 < \psi^{\prime\prime} (Y_s(t))\Phi_\xi (x,\lambda
\psi^\prime (Y_s (t) )), y_s-x  > dt,
\label{eq:Om1a}\end{equation}
\begin{equation}
\Omega_2 (x,\lambda)= -\lambda^{2-2r} \int_0^1 t\ < \psi^{\prime\prime} (Y_s(t))
(y_s-x), y_s-x  > dt.
\label{eq:Om2a}\end{equation}
\end{lemma}
  {\it Proof.} --
Let $\Omega_1(x,\lambda;t)$ be defined by formula (\ref{eq:Om1}) with $y_s$ replaced by $Y_s(t)$.
Since $\partial_t Y_s(t)=y_s-x$, we have that 
\[
\partial_t \Omega_1 (x,\lambda;t)=\lambda^{2-2r} < \psi^{\prime\prime} (Y_s(t))\Phi_\xi
(x,\lambda
\psi^\prime (Y_s (t) )), y_s-x  >.
\]
Integrating this equality and
taking into account that $\Omega_1 (x,\lambda;1)=\Omega_1 (x,\lambda)$ and $\Omega_1
(x,\lambda;0)\\=0$, we obtain (\ref{eq:Om1a}). Equality (\ref{eq:Om2a}) can be derived quite
similarly.$\quad\Box$

\noindent{\bf Remark}$\;$
If $\psi(x)$ is a linear function, then $\Omega(x,\lambda)=0$.

\begin{proposition}\label{S3} 
Function $(\ref{eq:GG})$ admits representation $(\ref{eq:S16})$, where  $\Omega (x,\lambda)$ as well
as all its derivatives $\Omega^{(\alpha)} (x,\lambda)=\partial^\alpha_x\Omega
(x,\lambda)$   are bounded uniformly in $\lambda$.
\end{proposition} 
{\it Proof.} -- 
By virtue of (\ref{eq:S8}) and the assumption $\Phi\in {\cal S}^r$, we have that $y_s-x$ and
$\Phi_\xi (x,\lambda
\psi^\prime (Y_s (t) ))$ are estimated by $C\lambda^{-1+r}$. Therefore functions (\ref{eq:Om1a}) and
(\ref{eq:Om2a}) are bounded uniformly in $\lambda$. 
 Differentiating equalities (\ref{eq:Om1a}) and (\ref{eq:Om2a}) and taking into account
(\ref{eq:S8df}), we see that the same is true for
all derivatives of the function
$\Omega(x,\lambda)$. $\quad\Box$

\begin{corollary}\label{S3C}
For all $\alpha$, the functions $\lambda^{-r} G^{(\alpha)} (x,\lambda)$  are bounded uniformly in
$\lambda$.
\end{corollary}

We shall need also a continuity of the  function $\Omega$ with respect to
variations of $\psi$.

\begin{proposition}\label{br1} 
Suppose that  functions $\psi(x)$ and $\breve{\psi}(x)$ satisfy Assumption~$\ref{psi}$. Set
\[
 \sigma_n= \max_{|\alpha|\leq n} \sup_{x\in U} 
 |\psi^{(\alpha)}(x )-\breve{\psi}^{(\alpha)}(x)|.
\]
 Let function   $\breve{\Omega}(x,\lambda)$ be
defined by equalities $(\ref{eq:Om1}) - (\ref{eq:Om})$  where $\breve{y}_s$ is the solution of
equation
\begin{equation}
 \breve{y}_s=x+ \Phi_\xi(x ,\lambda \breve{\psi}^\prime (\breve{y}_s)).
\label{eq:S5an}\end{equation}
 Then for any multi-index
$\alpha$ there exists a number $n(\alpha)$ such that, uniformly in $\lambda$, 
\[
 |\Omega^{(\alpha)}(x,\lambda)-\breve{\Omega}^{(\alpha)}(x,\lambda)|
\leq C_\alpha   \sigma_{n(\alpha)} ,\quad  x\in U.
\]
\end{proposition}

Let us start the proof with estimation of differences of the corresponding stationary points.

\begin{lemma}\label{br2}
For all $\alpha$  and sufficiently large $\lambda$,
\begin{equation}
 |y_s^{(\alpha)}(x,\lambda)-\breve{y}_s^{(\alpha)}(x,\lambda)|
\leq C_\alpha \lambda^{-1+r} 
\sigma_{|\alpha|+1}.
\label{eq:br3}\end{equation} 
\end{lemma}
{\it Proof.} --
Comparing equations (\ref{eq:S5a}) and (\ref{eq:S5an}), we find that
\begin{eqnarray*}
 |y_s-\breve{y}_s| =|\Phi_\xi(x ,\lambda \psi^\prime (y_s))-\Phi_\xi(x ,\lambda
\breve{\psi}^\prime (\breve{y}_s))|
\\
\leq C \lambda^{-1+r}| \psi^\prime (y_s) - \breve{\psi}^\prime (\breve{y}_s)|
\leq C_1 \lambda^{-1+r}\Bigl( |y_s-\breve{y}_s| + \sup_{x\in U} 
 |\psi^\prime(x )-\breve{\psi}^\prime(x)|\Bigr).
\end{eqnarray*}
This ensures estimate (\ref{eq:br3}) for $\alpha=0$. Differentiating (\ref{eq:S5a}),
(\ref{eq:S5an}) and using estimates (\ref{eq:S8df}) on
$y_s^{(\alpha)}(x,\lambda)$ and $\breve{y}_s^{(\alpha)}(x,\lambda)$, we can derive (\ref{eq:br3})
inductively for all $\alpha$.
$\quad\Box$

It follows from (\ref{eq:br3}) that for  any $\alpha$
\begin{equation}
| \psi^{(\alpha)}(Y_s)
-\breve{\psi}^{(\alpha)}(\breve{Y}_s)|\leq
\sup_{x\in U} 
 |\psi^{(\alpha)}(x )-\breve{\psi}^{(\alpha)}(x)| + C_\alpha|y_s-\breve{y}_s|\leq
\tilde{C}_\alpha(\sigma_{|\alpha|}+\sigma_1 \lambda^{-1+r}).
\label{eq:br4}\end{equation}
According to (\ref{eq:Om}) and (\ref{eq:Om1a}), (\ref{eq:Om2a}), to estimate the difference
$\Omega  -\breve{\Omega} $,
it suffices to use inequalities  (\ref{eq:br3}) for $\alpha=0$, (\ref{eq:br4}) for $|\alpha|=1,2$
and
\[
|\Phi_\xi (x,\lambda \psi^{\prime}(Y_s) )-\Phi_\xi (x,\lambda \breve{\psi}^{\prime}(\breve{Y}_s)
) | \leq C \lambda^{-1+r} | \psi^{\prime}(Y_s) -\breve{\psi}^{\prime}(\breve{Y}_s)|.
\]
Differentiating (\ref{eq:Om1a}), (\ref{eq:Om2a}) we can estimate  
$\Omega^{(\alpha)}  -\breve{\Omega}^{(\alpha)} $ quite in the same way.
This concludes the proof of Proposition~\ref{br1}.

{\bf 3.3.}
Recall that the phase function $G(x,\lambda)$ in (\ref{eq:I3A}) is defined by equations (\ref{eq:S5})
and equalities (\ref{eq:S3}), (\ref{eq:GG}). However
since we are interested only in the leading term of the asymptotics of the function
$ (Au_{\lambda,\varepsilon})(x )$ 
as $\lambda\rightarrow\infty$, we can neglect   the part of  
$G(x,\lambda)$ which tends to zero as $\lambda\rightarrow\infty$. This allows us to obtain a more
explicit expression for it.
 Thus,  in the case $r<1/2$, it follows from equality (\ref{eq:S16})
and Proposition~\ref{S3} that $G(x,\lambda)$ can be replaced  by   $\Phi (x,\lambda
\psi^\prime (x ))$.

In the case $r\in [1/2,2/3)$, we should keep the leading term of the asymptotics of the function
$\Omega(x,\lambda)$. To find it, we use relations (\ref{eq:S5a}), (\ref{eq:S8}) and, making an error
of order $O(\lambda^{-1+r})$, replace in (\ref{eq:Om1a}), (\ref{eq:Om2a}) the functions $y_s-x$ and 
$\Phi_\xi (x,\lambda \psi^\prime (Y_s(t)))$ by $\Phi_\xi (x,\lambda \psi^\prime (x ))$ and the
function $\psi^{\prime\prime} (Y_s(t))$ by $\psi^{\prime\prime} (x)$. It follows that
\begin{equation}
\Omega(x,\lambda)= 2^{-1} \lambda^{2-2r} < \psi^{\prime\prime}(x)
\Phi_\xi (x,\lambda \psi^\prime (x )), \Phi_\xi (x,\lambda \psi^\prime (x ))>+ O(\lambda^{-1+r}).
\label{eq:OM}\end{equation}
Inserting this expression into (\ref{eq:S16}), we obtain the two terms of the asymptotics of 
$G(x,\lambda)$ with an error  of order $O(\lambda^{-2+3r})$, which tends to zero if $r<2/3$ and
hence is negligible. We should keep, of course, more terms in the asymptotics of 
$G(x,\lambda)$ as $r$ increases.

 \section {the essential spectrum}

Our study of the essential spectrum of PDO with oscillating symbols (or amplitudes) relies on a 
construction of Weyl (singular) sequences. We seek these sequences in the form (\ref{eq:I2d}) and
proceed from Proposition~\ref{S5}. So our goal is to replace $G(x,\lambda)$ in (\ref{eq:I3A}) by
$G(x_0,\lambda)$. This requires a special choice of the function $\psi(x)=\psi(x,\lambda)$ which
will depend on the parameter $\lambda$.

{\bf 4.1.} 
We  start with an auxiliary construction which is non-trivial for $d>1$ only. 
Let ${\Bbb M}_n^{(d)}$ be the space of sequences $\Psi_n=\{\psi_\alpha
\},\;\psi_\alpha=\bar{\psi_\alpha}$, parametrized by multi-indices
$\alpha=(i_1,\ldots,i_n)$, $ 1\leq i_k\leq d,$ and symmetric with respect to all permutations of
indices $i_1,\ldots,i_n$. Thus, an element $\Psi_n$ is determined by numbers $\psi_\alpha$ for
$\alpha=(i_1,\ldots,i_n)$ with $1\leq
i_1\leq\ldots\leq i_n\leq d$. In particular, the set
${\Bbb M}_2^{(d)}$ can be identified with symmetric
$d\times d$- matrices. It is easy to check that the dimension $m_n^{(d)}$ of the space ${\Bbb
M}_n^{(d)}$ equals
\[
m_n^{(d)}=(n+d-1)! (n! (d-1)!)^{-1}.
\]
Below we  omit the index $``d"$.

 Assume that real numbers $t_1,\ldots,t_d$, such that
\[
||t||^2= \sum_{k=1}^d t_k^2\neq 0,
\]
and an element $F_{n-1}=\{f_\beta\}\in {\Bbb M}_{n-1}$ are given. Our goal here is to construct a
solution
$\Psi_n=\{\psi_\alpha\} \in {\Bbb M}_n$ of the system
\begin{equation}
 \sum_{k=1}^d \psi_{\beta,k} t_k =f_\beta,\quad \beta=(i_1,\ldots, i_{n-1}).
\label{eq:C2}\end{equation}
Of course,  (\ref{eq:C2}) contains $m_{n-1}$  equations for $m_n$ numbers $\psi_\alpha$ and 
$m_{n-1}<m_n$.  It is convenient to introduce the (annihilation) operator
 ${\bf T}_n:{\Bbb M}_n \rightarrow {\Bbb M}_{n-1}$
by the equality
\begin{equation}
({\bf T}_n \Psi_n)_\beta= \sum_{k=1}^d \psi_{\beta,k} t_k , 
\label{eq:AC1}\end{equation}
so that (\ref{eq:C2}) reads as ${\bf T}_n \Psi_n=F_{n-1}$.

Let us define also the (creation) operator
 ${\bf S}_n:{\Bbb M}_n \rightarrow {\Bbb M}_{n+1}$
by the equality
\[
({\bf S}_n \Psi_n)_{i_1,\ldots, i_{n+1}}=  \psi_{i_1,\ldots, i_n} t_{i_{n+1}} +
\psi_{i_1,\ldots, i_{n-1},i_{n+1}} t_{i_{n}}+\ldots+ \psi_{i_{n+1},i_2,\ldots, i_n} t_{i_1}.
\]
Clearly, 
\[
|| {\bf T}_n|| \leq C_n || t||,\quad || {\bf S}_n|| \leq C_n || t||.
\]
An easy computation shows that
\begin{equation}
 {\bf T}_{n+1}{\bf S}_n =|| t||^2 I_n+ {\bf S}_{n-1}  {\bf T}_n,
\label{eq:AC3}\end{equation}
where $I_n$ is the identity operator in ${\Bbb M}_n$. The relation between the creation and
annihilation operators becomes even more simple if one introduces the Fock space
\[
{\Bbb M}=\bigoplus_{n=0}^\infty {\Bbb M}_n, \quad {\Bbb M}_0={\Bbb R},
\]
and set
\[ 
{\bf T}=\bigoplus_{n=0}^\infty {\bf T}_n, 
\quad {\bf T}_0=\{0\}, \quad {\bf S}=\bigoplus_{n=0}^\infty {\bf S}_n.
\]
Then (\ref{eq:AC3}) for all $n$ are equivalent to the relation
\begin{equation}
 {\bf T}{\bf S} =|| t||^2 I+ {\bf S}  {\bf T}.
\label{eq:AC4}\end{equation}
Using (\ref{eq:AC4}) we can easily solve system (\ref{eq:C2}).

\begin{proposition}\label{C2}
For any $F_{n-1} \in {\Bbb M}_{n-1}$, $n\geq 1$, a solution of the equation ${\bf T}_n
\Psi_n=F_{n-1}$ can be constructed by the formula
\begin{eqnarray}
\Psi_n= || t||^{-2}\Bigl( {\bf S}- (2!)^{-1}|| t||^{-2} {\bf S}^2 {\bf T}
+ (3!)^{-1}|| t||^{-4} {\bf S}^3 {\bf T}^2-
\nonumber\\
\ldots -(-1)^n (n!)^{-1} || t||^{-2n+2} {\bf S}^n {\bf
T}^{n-1}\Bigr)F_{n-1}=: {\bf R}_{n-1} F_{n-1}.
\label{eq:AC5}\end{eqnarray}
In particular,
\begin{equation}
 ||{\bf R}_{n-1} ||\leq C_n  ||t||^{-1}.
\label{eq:C7}\end{equation}
\end{proposition}
{\it Proof.} --
It follows from (\ref{eq:AC4}) that
\[
 {\bf T}{\bf S}^p = p|| t||^2 {\bf S}^{p-1} + {\bf S}^p  {\bf T}
\]
so that
\[
 (p! || t||^{2p} )^{-1}{\bf T}{\bf S}^p {\bf T}^{p-1}F_{n-1} 
= 
((p-1)!  || t||^{2p-2})^{-1} {\bf S}^{p-1} {\bf T}^{p-1}F_{n-1}
+
 (p!  || t||^{2p})^{-1} {\bf S}^p {\bf T}^p F_{n-1}.
\]
Therefore applying the operator ${\bf T}$ to equality (\ref{eq:AC5}) we see that all terms except
the first, which is $F_{n-1}$, and the last, which is
\[
-(-1)^n (n!  || t||^{2n})^{-1} {\bf S}^n {\bf T}^n F_{n-1}=0,
\]
cancel each other. Hence ${\bf T}\Psi_n=F_{n-1}.\quad\Box$

{\bf 4.2.}
Let the phase function $\Phi\in{\cal S}^r$ be given. For a function $\psi(x,\lambda)$ depending on
the parameter $\lambda$, define the function
$G(x,\lambda)$  by formulas (\ref{eq:S3}), (\ref{eq:S5}) (where $\psi=\psi(x,\lambda)$) and 
(\ref{eq:GG}). This definition is correct for sufficiently large $\lambda$ as long as the family
$\psi(\cdot,\lambda)$ satisfies Assumption~\ref{psi}.

Choose a point $x_0\in X$ and some $n=1,2,\ldots $. Our goal is to find functions
$\psi(x,\lambda)$ such that
\begin{equation}
(\partial^\alpha G)(x_0,\lambda)=0,\quad |\alpha|=1,\ldots , n,\quad \partial=\partial_x.
\label{eq:G1}\end{equation}
We seek $\psi(x,\lambda)$ as a polynomial of degree $n+1$: 
\begin{equation}
\psi(x,\lambda)=\sum_{1\leq |\alpha|\leq n+1}(\alpha !)^{-1} \psi_\alpha(\lambda)(x-x_0)^\alpha,
\label{eq:G2}\end{equation}
 where, of course, 
$\psi_\alpha(\lambda)=(\partial^\alpha \psi)(x_0,\lambda)$. Below we fix   a vector
$\xi_0=\psi^\prime(x_0)\neq 0$, which does not depend on $\lambda$. Let us denote by
$\Psi_k=\Psi_k(\lambda),\; k=2,\ldots,n+1$, the collection of
$\{\psi_\alpha(x_0,\lambda)\}$ for all $|\alpha|=k$.

 By virtue of Proposition~\ref{S3}, the asymptotics of $(\partial^\alpha G)(x,\lambda)$ for large
$\lambda$ is determined  by the term $\partial^\alpha \Phi(x,\lambda\psi^\prime(x,\lambda ))$.
We need to single out the derivative of the highest order (which equals to $|\alpha|+1$) of this
function. Denote by $\Phi^{(\alpha)}_x(x,\xi)$ the derivative of $\Phi (x,\xi)$  in the variable
$x$ of order $\alpha$.

\begin{lemma}\label{G1}
Let $\Phi \in{\cal S}^r$. Then for all $\alpha$
\begin{equation}
\partial^\alpha \Phi(x,\lambda\psi^\prime(x,\lambda ))=\lambda\sum_{i=1}^d
\Phi_{\xi_i}(x,\lambda\psi^\prime(x,\lambda)) (\partial_i \partial^\alpha \psi)(x )+
\lambda^r F_\alpha(x,\lambda ),
\label{eq:G3}\end{equation} 
where 
\begin{eqnarray} 
F_\alpha(x,\lambda )= \lambda^{-r}\Phi^{(\alpha)}_x (x,\lambda\psi^\prime(x,\lambda ))
\nonumber\\
 +\sum_
{\scriptstyle 2\leq |\beta_l|\leq |\alpha|, 1\leq l \leq|\alpha|}
 \lambda^{l-r} f_{\beta_1,\ldots,\beta_l}(x,\lambda\psi^\prime(x,\lambda ))
\psi^{(\beta_1)}(x,\lambda )\ldots \psi^{(\beta_l)}(x,\lambda )
\label{eq:G4}\end{eqnarray}
with functions $f_{\beta_1,\ldots,\beta_l} \in{\cal S}^{r-l}$.
\end{lemma}

 A proof  can be  easily obtained by induction in the order of derivatives.

The following assertion is a direct consequence of (\ref{eq:G4}). 

\begin{lemma}\label{GH2}
If $\psi(x,\lambda)$ is defined by $(\ref{eq:G2})$, then
\[
\lambda^{-r}F_\alpha(x_0,\lambda )=P_\alpha(\Psi_2(\lambda) ,\ldots,\Psi_{|\alpha|}(\lambda)
;\lambda),
\]
where $P_\alpha(\Psi_2  ,\ldots,\Psi_{|\alpha|};\lambda) $ is a polynomial of $\Psi_2
,\ldots,\Psi_{|\alpha|}$ with uniformly in $\lambda$ bounded coefficients.
\end{lemma}

Recall that the function 
 $\Omega(x,\lambda) $ is defined by equalities (\ref{eq:Om1}) - (\ref{eq:Om}). For any $\lambda$,  
the derivative $ \Omega^{(\alpha)}(x_0,\lambda)$  is
 determined by sequences $\Psi_k=\Psi_k(\lambda)$, that is
\[
\Omega^{(\alpha)}(x_0,\lambda)=\Omega_\alpha(\Psi_2(\lambda),\ldots,
\Psi_{n+1}(\lambda);\lambda)
\]
for some function $\Omega_\alpha$. The following result is a direct consequence of 
Propositions~\ref{S3} and \ref{br1}.

\begin{lemma}\label{GH1}
The functions $\Omega_\alpha(\Psi_2 ,\ldots,
\Psi_{n+1} ;\lambda)$ satisfy the estimates
\begin{equation} 
|\Omega_\alpha(\Psi_2 ,\ldots,
\Psi_{n+1} ;\lambda)|\leq C 
\label{eq:XY2}\end{equation}
\begin{equation}
 |\Omega_\alpha(\Psi_2^\prime ,\ldots,
\Psi_{n+1}^\prime ;\lambda) - \Omega_\alpha(\Psi_2^{\prime\prime} ,\ldots,
\Psi_{n+1}^{\prime\prime} ;\lambda)|\leq C
\sum_{k=2}^{n+1}|\Psi_k^\prime -\Psi_k^{\prime\prime}|
\label{eq:XY3}\end{equation}
for $\Psi_2 ,\ldots,\Psi_{n+1} $ from any compact subsets and sufficiently large $\lambda$.
\end{lemma}

We say that a function $\Omega_\alpha$ obeying (\ref{eq:XY2}), (\ref{eq:XY3}) satisfies the Lipschitz
condition in $\Psi_2 ,\ldots,\Psi_{n+1} $ uniformly in $\lambda$.

Comparing  (\ref{eq:S16}) and (\ref{eq:G3}) and
putting $x=x_0$ we can now rewrite equations (\ref{eq:G1}) as 
\begin{equation}
\sum_{i=1}^d \psi_{\alpha,i}(\lambda) t_i(\lambda) +  
P_\alpha(\Psi_2(\lambda),\ldots,\Psi_{|\alpha|}(\lambda);\lambda) + \lambda^{-1+r}
\Omega_\alpha(\Psi_2(\lambda),\ldots,\Psi_{n+1}(\lambda);\lambda) =0,  
\label{eq:G5}\end{equation}
where $1\leq|\alpha|\leq n$ and
\begin{equation}
t_i(\lambda)= \lambda^{1-r}\Phi_{\xi_i}(x_0,\lambda \xi_0).
\label{eq:G6}\end{equation}
Below we often omit in notation the dependence of different  functions on the parameter $\lambda$
which is supposed to be large.
 We treat (\ref{eq:G5}) as a system of equations for ``vectors"
$\Psi_2,\ldots,\Psi_{n+1}$. Using definition (\ref{eq:AC1}) of
the operator
${\bf T}_{k+1}={\bf T}_{k+1}(\lambda)$ we write the set of equations (\ref{eq:G5}) with $|\alpha|=k$
in the vector notation
\begin{equation}
{\bf T}_{k+1}(\lambda)\Psi_{k+1}  + P_k(\Psi_2 ,\ldots,\Psi_k ;\lambda) +
\lambda^{-1+r}
\Omega_k(\Psi_2 ,\ldots,\Psi_{n+1} ;\lambda) =0,\quad k=1,\ldots, n.  
\label{eq:G5v}\end{equation}
We emphasize that $P_1(\lambda)=\lambda^{-r}\Phi_x(x_0,\lambda\xi_0)$  does not
depend on sequences
$\Psi_2,\ldots,\Psi_{n+1}$. Thus, we have the following

\begin{lemma}\label{GH3}
System of equations $(\ref{eq:G1})$ for polynomial $(\ref{eq:G2})$ is equivalent to system
$(\ref{eq:G5v})$ for vectors $\Psi_2,\ldots,\Psi_{n+1}$.
\end{lemma}

 Suppose that $x_0$ and $\xi_0=\psi^\prime(x_0)$ are chosen in such a way
that 
\begin{equation}
|\Phi_\xi(x_0,\lambda \xi_0) | \geq c \lambda^{r-1}.
\label{eq:G7}\end{equation}
Then the vector $t(\lambda)$ with components (\ref{eq:G6}) satisfies $||t(\lambda)||\geq c>0$.
Let the operator ${\bf R}_k={\bf R}_k(\lambda)$ be defined by (\ref{eq:AC5}).
Set $\Psi_{k+1}={\bf R}_k \tilde{\Psi}_k$ and
\[
 \tilde{P}_k (\tilde{\Psi}_1,\ldots,\tilde{\Psi}_{k-1};\lambda) = -P_k ({\bf
R}_1\tilde{\Psi}_1,\ldots,{\bf R}_{k-1}\tilde{\Psi}_{k-1};\lambda),
\]
\[
\tilde{\Omega}_k (\tilde{\Psi}_1,\ldots,\tilde{\Psi}_n;\lambda)=-  \Omega_k
 ({\bf R}_1\tilde{\Psi}_1,\ldots,{\bf R}_n\tilde{\Psi}_n;\lambda).  
\]
 By Proposition~\ref{C2},
\[
{\bf T}_{k+1} {\bf R}_k \tilde{\Psi}_k =\tilde{\Psi}_k 
\]
so that system (\ref{eq:G5v}) can be rewritten as
\begin{equation}
 \tilde{\Psi}_k =\tilde{P}_k (\tilde{\Psi}_1,\ldots,\tilde{\Psi}_{k-1};\lambda) +
\lambda^{-1+r}\tilde{\Omega}_k
(\tilde{\Psi}_1,\ldots,\tilde{\Psi}_n;\lambda),\quad k=1, \ldots, n,  
\label{eq:G8}\end{equation}
for $\tilde{\Psi}_1,\ldots,\tilde{\Psi}_n$ who will depend of course  on $\lambda$.
Taking also into account estimate    (\ref{eq:C7}) and Lemmas~\ref{GH1}, \ref{GH2}, we can
formulate an intermediary result.

\begin{lemma}\label{GH4}
Functions $ \tilde{\Psi}_k$ and $\tilde{\Omega}_k$ satisfy the Lipschitz
condition in the variables
$\tilde{\Psi}_1,\ldots,\newline\tilde{\Psi}_{k-1}$ and $\tilde{\Psi}_1,\ldots,\tilde{\Psi}_n$,
respectively, uniformly in $\lambda$. If system $ (\ref{eq:G8})$ is
fulfilled for $\tilde{\Psi}_1,\ldots,\tilde{\Psi}_n$, then  system $ (\ref{eq:G5v})$ is fulfilled for
$\Psi_2={\bf R}_1\tilde{\Psi}_1,\ldots,\Psi_{n+1} = {\bf R}_n\tilde{\Psi}_n$.
\end{lemma}

Remark that, up to a small term $\lambda^{-1+r}\tilde{\Omega}_k$, system
(\ref{eq:G8}) has a triangular structure, i.e. $\tilde{P}_k $ depends on
$\tilde{\Psi}_1,\ldots,\tilde{\Psi}_{k-1}$ only. This allows us to solve system (\ref{eq:G8}) by
iterations starting from
\[
 \tilde{\Psi}_1^{(0)}=  \tilde{P}_1 (\lambda),\quad 
\tilde{\Psi}_k^{(0)} =\tilde{P}_k
(\tilde{\Psi}_1^{(0)},\ldots,\tilde{\Psi}_{k-1}^{(0)};\lambda), \quad  k =2,\ldots,n.
\]
Set
\begin{equation}
 \tilde{\Psi}_1^{(p+1)} =\tilde{P}_1 (\lambda) + \lambda^{-1+r}
\tilde{\Omega}_1 (\tilde{\Psi}_1^{(p)},\ldots,\tilde{\Psi}_n^{(p)};\lambda),
\label{eq:G12}\end{equation} 
\begin{equation}
\tilde{\Psi}_k^{(p+1)} =\tilde{P}_k
(\tilde{\Psi}_1^{(p+1)},\ldots,\tilde{\Psi}_{k-1}^{(p+1)};\lambda) + \lambda^{-1+r}
\tilde{\Omega}_k (\tilde{\Psi}_1^{(p)},\ldots,\tilde{\Psi}_n^{(p)};\lambda), \quad  k =2,\ldots,n.
\label{eq:G13}\end{equation}

\begin{lemma}\label{G2}
 For all $p\geq 1$ and $1\leq k\leq n$
\begin{equation}
|\tilde{\Psi}_k^{(p)}-\tilde{\Psi}_k^{(p-1)}|\leq C \lambda^{-p(1-r)}.
\label{eq:G14}\end{equation} 
\end{lemma} 
{\it Proof.}--
Suppose that (\ref{eq:G14}) holds for some $p=p_0$. Let us check it for $p=p_0+1$. Remark, first,
that according to  equality (\ref{eq:G12})
\[
 \tilde{\Psi}_1^{(p+1)}-\tilde{\Psi}_1^{(p)}= \lambda^{-1+r}
\Bigl(\tilde{\Omega}_1 (\tilde{\Psi}_1^{(p)},\ldots,\tilde{\Psi}_n^{(p)};\lambda)
-\tilde{\Omega}_1 (\tilde{\Psi}_1^{(p-1)},\ldots,\tilde{\Psi}_n^{(p-1)};\lambda)\Bigr).
\]
By Lemma~\ref{GH4}, the function $\tilde{\Omega}_1$ is uniformly in $\lambda$ Lipschitz
 continuous so that the right-hand side here is bounded by 
\begin{equation}
C \lambda^{-1 +r} \sum_{1\leq k \leq n} |\tilde{\Psi}_k^{(p)}-\tilde{\Psi}_k^{(p-1)}|.
\label{eq:G15}\end{equation}
This does not exceed $C \lambda^{-(p+1)(1-r)}$ by our inductive assumption.
 Supposing, further, that
\begin{equation}
 |\tilde{\Psi}_k^{(p+1)}-\tilde{\Psi}_k^{(p)}|\leq C \lambda^{-(p+1)(1-r)},
\label{eq:G16}\end{equation}
we check it for $k=k_0+1$. According to  equality (\ref{eq:G13}),
\begin{eqnarray}
 \tilde{\Psi}_k^{(p+1)}-\tilde{\Psi}_k^{(p)}= 
\Bigl(\tilde{P}_k(\tilde{\Psi}_1^{(p+1)},\ldots,\tilde{\Psi}_{k-1}^{(p+1)},\lambda) -
\tilde{P}_k(\tilde{\Psi}_1^{(p)},\ldots, \tilde{\Psi}_{k-1}^{(p)},\lambda)\Bigr)
\nonumber\\
+\lambda^{-1+r} \Bigl(\tilde{\Omega}_k(\tilde{\Psi}_1^{(p)},\ldots, \tilde{\Psi}_n^{(p)},\lambda) -
\tilde{\Omega}_k( \tilde{\Psi}_1^{(p-1)},\ldots, \tilde{\Psi}_n^{(p-1)},\lambda)\Bigr).
\label{eq:G17}\end{eqnarray}
Since, again by Lemma~\ref{GH4}, the function $\tilde{P}_k$ is uniformly in $\lambda$ Lipschitz
 continuous, the first term in the right-hand side is estimated by 
\[
C  \sum_{l=1}^{k-1} | \tilde{\Psi}_l^{(p+1)}- \tilde{\Psi}_l^{(p)}|\leq C_1 \lambda^{-(p+1)(1 -r)}.
\]
Here we used conjecture (\ref{eq:G16}). The second term in the right-hand side of (\ref{eq:G17})
is estimated, as before, by (\ref{eq:G15}). This proves (\ref{eq:G16}) for $k=k_0+1$ and,
consequently, (\ref{eq:G14}) for $p=p_0+1$. The same arguments show, of course, also that 
(\ref{eq:G14}) is fulfilled for $p=1.\quad\Box$

It follows from  (\ref{eq:G14}) that sequences $\tilde{\Psi}_k^{(p)},\: 1\leq k\leq n$,
have finite limits $\tilde{\Psi}_k$ as $p\rightarrow\infty$. Passing to these limits in relations 
(\ref{eq:G13}) we see that $\tilde{\Psi}_1,\ldots,\tilde{\Psi}_n$ satisfy system (\ref{eq:G8}). By
Lemma~\ref{GH4}, this implies that (\ref{eq:G5v}) is satisfied for
$\Psi_2={\bf R}_1\tilde{\Psi}_1,\ldots,\Psi_{n+1}={\bf R}_n\tilde{\Psi}_n$. Thus, we arrive at the
following assertion.

\begin{proposition}\label{G3}
Suppose that condition $(\ref{eq:G7})$ is satisfied. Then for any $n\geq 2$ and sufficiently large
$\lambda$, a solution $\Psi_2(\lambda),\ldots, \Psi_{n+1}(\lambda)$ of system
$(\ref{eq:G5})$ exists. For the corresponding function $(\ref{eq:G2})$ and function $G(x,\lambda)$
defined by $(\ref{eq:S3})$, $(\ref{eq:S5})$ and
$(\ref{eq:GG})$, equations $(\ref{eq:G1})$ are fulfilled.
\end{proposition}

{\bf 4.3.}
Suppose that equality (\ref{eq:BC2}) is satisfied with $b\in {\cal S}^0$.
Let functions $u_{\lambda,\varepsilon}$ be defined by equality (\ref{eq:I2d}).
According to Proposition~\ref{S5}, the
asymptotics of $Au_{\lambda,\varepsilon}$ as   $ \lambda \rightarrow\infty$,  
 $\varepsilon\lambda^{1-r}\rightarrow\infty$ is given by relation (\ref{eq:I3A}).
Let $\psi(x,\lambda)$ be polynomial (\ref{eq:G2}) with coefficients $\psi_\alpha(\lambda)$
satisfying (\ref{eq:G5})  for $1\leq |\alpha|\leq n$ so that equalities (\ref{eq:G1}) are fulfilled.
Then it follows from Corollary~\ref{S3C} that
\begin{equation}
|G(x,\lambda) - G(x_0,\lambda)| \leq C \lambda^r |x-x_0|^{n+1}.
\label{eq:G18}\end{equation}
The right-hand side here tends to zero if $\lambda^r \varepsilon^{n+1}\rightarrow 0$. This is
compatible with the condition $\lambda^{-1+r} \varepsilon^{-1}\rightarrow 0$ if
 \begin{equation}
n+1> r(1-r)^{-1}.
\label{eq:nr}\end{equation}
Combining (\ref{eq:I3A}) with (\ref{eq:I3An}) and (\ref{eq:G18}), we obtain

\begin{proposition}\label{G3h} 
Let functions $u_{\lambda,\varepsilon}$ be defined by relations $(\ref{eq:I2d})$, $(\ref{eq:G2})$ and
$(\ref{eq:G5})$. Then
\begin{equation}
  ||Au_{\lambda,\varepsilon}- e^{iG(x_0,\lambda)} b(x_0, \lambda
\xi_0 )u_{\lambda,\varepsilon}||
\rightarrow 0
\label{eq:G19}\end{equation}
 as $\lambda\rightarrow \infty$ and
$\lambda^{-1+r}\varepsilon^{-1}\rightarrow 0,
\lambda^r\varepsilon^{n+1}\rightarrow 0$.
\end{proposition}

Suppose additionally that
\begin{equation} 
|\Phi (x_0,\lambda \xi_0)|\geq c \lambda^r,\quad c>0.
\label{eq:G20}\end{equation}
Equality (\ref{eq:S16}) and Proposition~\ref{S3} imply that the same inequality is true
for $G(x_0,\lambda)$ and, in particular,
\[ 
 \lim_{\lambda\rightarrow\infty}G(x_0,\lambda)=\infty\quad {\rm or}\quad
\lim_{\lambda\rightarrow\infty}G(x_0,\lambda)=-\infty.
\]
Since $G(x_0,\lambda)$ is a continuous function of $\lambda$, for any $\mu_1=e^{i\theta}$, we can
find a sequence $\lambda_p\rightarrow\infty$ such that $G(x_0,\lambda_p)=\theta +2\pi p$ or
$G(x_0,\lambda_p)=\theta -2\pi p$. Then 
$e^{iG(x_0,\lambda_p)}=\mu_1$. Under assumption (\ref{eq:nr}) we can set 
\[
\varepsilon_p=\lambda_p^{-s}\quad {\rm for} \quad (n+1)^{-1} r< s < 1-r
\]
so that
$
\lambda_p^{-1+r}\varepsilon_p^{-1}\rightarrow 0,\;
\lambda_p^r \varepsilon_p^{n+1}\rightarrow 0
$
as $p\rightarrow\infty$. 
If condition (\ref{eq:I1b}) holds, then
it follows from (\ref{eq:G19}) that for functions $u_p=u_{\lambda_p,\varepsilon_p}$
\[
\lim_{p\rightarrow\infty} || Au_p-\mu u_p||=0,\quad \mu=\mu_1\mu_0,
\]
and hence $\mu$ belongs to the spectrum of $A$.  Let us formulate the result obtained.

\begin{theorem}\label{G4}
Let the symbol $a(x,\xi)$ of a PDO $(\ref{eq:BC1})$ be compactly supported in $x$ and satisfy
conditions $(\ref{eq:BC2})$ with $m=0$.
 Suppose that for some point $x_0\in X, \xi_0\neq 0$ conditions $(\ref{eq:I1b}),
(\ref{eq:G7})$ and $(\ref{eq:G20})$
are satisfied. Then  the spectrum of the operator $A $ in the space $L_2(X)$ covers the
circle  ${\Bbb T}_\kappa$, where $\kappa=|\mu_0|$.
\end{theorem}

Combining this result with Theorem~\ref{SP4}, we can generalize Theorem~\ref{G4}
 to PDO defined by oscillating amplitudes. 

\begin{theorem}\label{G4A}
 Let the amplitude $a(x,x^\prime,\xi)$ of a PDO $(\ref{eq:RH1})$ be compactly supported in $x$ and
$x^\prime$  and satisfy conditions $(\ref{eq:I1})$ with $m=0$. Suppose that
\[
\lim_{\lambda\rightarrow\infty}{\bf b}(x_0,x_0,\lambda\xi_0)=\mu_0 \neq 0, \quad
|\Theta (x_0,x_0,\lambda\xi_0)| \geq c \lambda^r, \quad
|\Theta_\xi (x_0,x_0,\lambda\xi_0)| \geq c \lambda^{r-1} 
\]
 for some point $x_0\in X, \xi_0\neq 0$ and $ c>0$. Then  the spectrum of the
operator $A$ in the space $L_2(X)$ covers the circle ${\Bbb T}_\kappa$, where $\kappa=|\mu_0|$.
\end{theorem}

Of course, assumptions of  Theorems~\ref{G4} and \ref{G4A} are fulfilled for
asymptotically homogeneous in $\xi$ functions $\Phi(x, \xi)$ and $\Theta(x,x^\prime,\xi)$.

\noindent{\bf Remark}$\;$
If $K$ is a compact operator, then, under the assumptions of Theorems~\ref{G4} and \ref{G4A},
the spectrum of the operator $A+K$   covers the circle ${\Bbb T}_\kappa$. This is a general result
on compact perturbations but it follows also from the proofs of these theorems  because the
constructed sequence $u_p$ converges weakly to zero as $p\rightarrow\infty$.

{\bf 4.4.}
Let us, finally, discuss particular cases $r<1/2$ and, more generally,  $r<2/3$. If $r<1/2$, then
(\ref{eq:nr}) holds for $n=0$ and $\psi(x)$ can be defined by formula (\ref{eq:1.9}). In this case
$G(x_0,\lambda)=\Phi(x_0,\lambda\xi_0)$ in (\ref{eq:G19}) and we
can  omit condition (\ref{eq:G7}) in Theorem~\ref{G4} (or a similar condition on $\Theta_\xi$ in
Theorem~\ref{G4A}).

If $r\in[1/2,2/3)$, then (\ref{eq:nr}) requires $n=1$ so that
\begin{equation} 
 \psi(x,\lambda)=<\xi_0,x-x_0> + 2^{-1} < \Psi_2(\lambda) (x-x_0), x-x_0>,
\label{eq:pc1}\end{equation}
where $\Psi_2(\lambda)$ is defined by equations (\ref{eq:G1}) for $n=1$. Actually,
$\Psi_2(\lambda)$ can be chosen as a simple but approximate solution of these equations.
 According to Proposition~\ref{S3},
\begin{equation}
 G(x,\lambda)= G(x_0,\lambda)+ <G^\prime(x_0,\lambda), x-x_0>+ O(\lambda^r |x-x_0|^2 ) 
\label{eq:pc2}\end{equation}
and 
\begin{equation}
G^\prime(x_0,\lambda)= G_1(\lambda)+  O(\lambda^{-1+2r}), \quad {\rm where}\quad
 G_1(\lambda)=\Phi_x (x_0,\lambda \xi_0)+ \lambda\Psi_2(\lambda)\Phi_\xi (x_0,\lambda \xi_0).
\label{eq:pc2x}\end{equation}
Suppose now that $ G_1(\lambda)=0$ and choose $\varepsilon=\lambda^{-s}$ with
\[
\max \{2r-1,r/2\}<s< 1-r,
\]
 which is possible if $r<2/3$. Then, by virtue of (\ref{eq:pc2}),
(\ref{eq:pc2x}), $G(x,\lambda)$ can be replaced by $G(x_0,\lambda)$ in (\ref{eq:I3A}) so that
(\ref{eq:G19}) holds. Thus, if $\Psi_2(\lambda)$ satisfies the equation $ G_1(\lambda)=0$, then
$\psi(x,\lambda)$ in (\ref{eq:I2d}) can be defined by equality (\ref{eq:pc1}). In this case
(\ref{eq:G19}) holds with $G(x_0,\lambda)$ given by formulas (\ref{eq:S16}), (\ref{eq:OM}) where
$x=x_0$.

\section {Integral kernels}

We treat here  PDO   as integral operators in one of curvilinear coordinates, whose kernels are PDO
in remaining coordinates. We distinguish a class of operators with continuous kernels so that, in
particular,  diagonal values of kernels are well-defined. In this section we consider PDO
$A:C_0^\infty(X)\rightarrow C^\infty(X)$  defined by formula (\ref{eq:RH1}). We suppose that the
amplitude
  ${\bf a}\in{\cal S}^{m}_{\rho, \delta,\delta}(X\times X \times {\Bbb R}^d)$ for some $m$ and
$\rho>0$,
$\delta<1$ but do not make any special assumptions of the type (\ref{eq:I1}).

{\bf 5.1.}
Below we need a formula of  change of variables for PDO defined by their amplitudes. For a
diffeomorphism
$\kappa: X\rightarrow Y$,  define the operator $F_\kappa$ by the relation
\begin{equation}
 (F_\kappa u)(y)= |\det \kappa^\prime(x)|^{-1/2} u(x),\quad {\rm where}\quad y =\kappa(x).
\label{eq:IK1}\end{equation} 
In view of our applications we introduced the factor $|\det \kappa^\prime(x)|^{-1/2}$ so that 
$F_\kappa:L_2(X)\rightarrow L_2(Y)$ is a unitary operator. 
Let $ G(x,x^\prime)$ be a $C^\infty$-operator-function satisfying
\begin{equation}
\kappa(x)-\kappa(x^\prime)=G(x,x^\prime)(x-x^\prime).
\label{eq:rh2s}\end{equation}
Then $ G(x,x)=\kappa^\prime (x)$ and $\det G(x,x^\prime)\neq 0$ in some neighbourhood $\Omega$ of
the diagonal $x=x^\prime$. One of  solutions of (\ref{eq:rh2s}) is given by the equality
\begin{equation}
 G(x,x^\prime)=\int_0^1 \kappa^\prime (x+t(x^\prime-x))dt.
\label{eq:rh2}\end{equation}

 Let  $\chi\in C_0^\infty (X\times X)$ be such that ${\rm supp}\: \chi\in \Omega$,
$\chi(x,x^\prime)=1$ in a neighbourhood of the diagonal and $\chi_0=1-\chi$. Then 
\begin{equation}
 (Au,v)_{L_2(X)}=(2\pi)^{-d}\int_{X} \int_{X} \int_{{\Bbb R}^d} e^{i<x-x^\prime,\xi>} {\bf
a}(x,x^\prime,\xi) u(x^\prime)\overline{v(x)} dx dx^\prime d\xi
\label{eq:RH1a}\end{equation}
 is a sum of the two integrals corresponding to $\chi_0{\bf a}$ and $\chi
{\bf a}$. The first term is 
$(K_0 u,v)_{L_2(X)}$, where $K_0$ is an integral operator with  $C^\infty$-kernel which equals to
zero in a neighbourhood of the diagonal. In the second integral we
 change  the variables 
\begin{equation}
x=\kappa^{-1}(y),\quad x^\prime=\kappa^{-1}(y^\prime), \quad \xi=\;^t G(x,x^\prime)\eta
\label{eq:rh3x}\end{equation}
  and set  $\tilde{u}=F_\kappa u,\; \tilde{v}=F_\kappa v $. According to (\ref{eq:rh2s}),
\begin{equation}
 <x-x^\prime,\xi> =< \kappa(x)-\kappa(x^\prime),\eta> 
\label{eq:rh3}\end{equation}
so that
\[
 (Au,v)_{L_2(X)}= (K_0u,v)_{L_2(X)}+(2\pi)^{-d}\int_{Y} \int_{Y} \int_{{\Bbb R}^d} e^{i<y-y^\prime,
\eta>} {\bf \tilde{a}}(y, y^\prime ,\eta) \tilde{u}(y^\prime)\overline{\tilde{v}(y)}dy dy^\prime
d\eta,
\]
 where
\begin{equation}
{\bf \tilde{a}}(y, y^\prime ,\eta) = \alpha(x,x^\prime) {\bf a} (x, x^\prime ,\xi)
\label{eq:rh6}\end{equation}
and
\begin{equation}
\alpha(x,x^\prime) =|\det \kappa^\prime(x)\; \det \kappa^\prime(x^\prime)|^{-1/2} \; |\det
G(x,x^\prime)| \chi(x,x^\prime).
\label{eq:rh7}\end{equation} 
Function (\ref{eq:rh6}) belongs to the class
${\cal S}^{m}_{\rho,\tilde{\delta},\tilde{\delta}}(Y\times Y\times{\Bbb R}^d)$ with
$\tilde{\delta}=\max\{\delta, 1-\rho\}$. Thus, we arrive at the following

\begin{proposition}\label{df}
Suppose that 
  ${\bf a}\in{\cal S}^{m}_{\rho, \delta,\delta}$, where
$\rho>0$, $\delta<1$  and $\rho+\delta\geq 1$. For a diffeomorphism $\kappa:X\rightarrow Y$, define
the operator   $F=F_\kappa$ by equality
$(\ref{eq:IK1})$ and let $G(x,x^\prime)$ satisfy equation
 $(\ref{eq:rh2s})$ $($for example, $G=G_\kappa$ can be defined by $(\ref{eq:rh2}))$. Let $\tilde{A}$
be the PDO with the amplitude
${\bf \tilde{a}}(y, y^\prime ,\eta)$
determined by equalities  $(\ref{eq:rh3x})$, $(\ref{eq:rh6})$ and $(\ref{eq:rh7})$
$($so that, in particular, ${\bf \tilde{a}}\in{\cal S}^{m}_{\rho,\delta,\delta})$. Then $FAF^{-1}
-\tilde{A}$ is an integral operator with
$C^\infty$-kernel which equals to zero in a neighbourhood of the diagonal $y=y^\prime$.
\end{proposition}

{\bf 5.2.}
Let us give an abstract definition of an integral operator  in the space $\mbox{\sf
H}=L_2(\Lambda;{\cal N})$ of  vector-functions $u(\lambda)$ defined on an interval
$\Lambda\subset{\Bbb R}$ and taking values in an auxiliary Hilbert space ${\cal N}$. Suppose that a 
set ${\cal D}\subset{\cal N}$ is  dense  in ${\cal N}$ and introduce the space ${\cal
G}=C_0^\infty(\Lambda;{\cal D})$ of infinitely differentiable compactly supported vector-functions
with values in ${\cal D}$. Clearly, ${\cal G}$ is dense in $\mbox{\sf H}$. Consider an operator
$A:{\cal G}\rightarrow {\cal G}^\prime$ where ${\cal G}^\prime$ is the dual space to ${\cal
G}$. Below we write $(w,v)_{\mbox{\sf H}}$ for elements $v\in{\cal G}$, $w\in{\cal G}^\prime$ so
that, strictly speaking, $(\cdot,\cdot)_{\mbox{\sf H}}$ is the duality symbol.

\begin{definition}\label{RS}
 Let  $A^\natural (\mu,\nu):{\cal D}\rightarrow {\cal D}^\prime$ be a
continuous operator-function of variables $\mu,\nu\in \Lambda$. We say that 
 $A^\natural (\mu,\nu)$ is  kernel of an operator $A $ if for any
$u,v\in {\cal G}$
\begin{equation}
 (Au,v)_{\mbox{\sf H}}= \int_{\Lambda} \int_{\Lambda}   (A^\natural (\mu,\nu) u(\nu), v(\mu
))_{\cal N} d\mu d\nu. 
\label{eq:RH1c}\end{equation} 
\end{definition}

 The bilinear form of the operator $A^\natural (\mu,\nu)$ can be constructed in terms of 
$A$. Indeed, 
let $g\in C^\infty ({\Bbb R})\cap L_1 ({\Bbb R})$,
$\int_{-\infty}^\infty g(t)dt =1$, $\psi_\mu\in C_0^\infty(\Lambda)$, $\psi_\mu(\lambda)=1$ in a
neighbourhood of the point $\mu$ and
\[
\psi_{\varepsilon,\mu}(\lambda)=\varepsilon^{-1}
g(\varepsilon^{-1}(\lambda-\mu))\psi_\mu(\lambda).
\]
 Then for any  $u, v\in{\cal D}$ the double limit exists 
\[
\lim_{\varepsilon,\eta\rightarrow 0} (A (\psi_{\varepsilon,\nu}u), \psi_{\eta,\mu}v)_{\mbox{\sf H}}=
(A^\natural (\mu,\nu) u, v)_{\cal N}.
\]

Thus kernel $A^\natural (\mu,\nu)$ of an operator $A$ is necessarily unique. Of course, only
operators from a rather restricted class may have continuous kernels. For example, kernel of a
Hilbert-Schmidt operator is not, in general continuous. On the other hand, an unbounded operator
(even defined as a mapping from ${\cal G}$ into ${\cal G}^\prime$ only) may have continuous kernel
in the sense of Definition~\ref{RS}. We point out that, for an operator with continuous kernel,
diagonal values 
$A^\natural (\lambda,\lambda):{\cal D}\rightarrow {\cal D}^\prime$ 
are well defined and continuous in $\lambda\in\Lambda$.

For a PDO $A:C_0^\infty(X)\rightarrow C^\infty(X)$, we consider its kernel in a direct integral
constructed with respect to the operator of multiplication by some smooth function $P(x)$.
Suppose first that the set $X\subset {\Bbb R}^d$ has a special structure. Let $x=(x_0,x_d)$
where $x_0=(x_1,\ldots, x_{d-1})$, and define the mapping
$\kappa: X \rightarrow Y $ by the equalities
\begin{equation} 
 y_0=x_0,\quad y_d=P(x).
\label{eq:AX1}\end{equation} 
Assume that $\kappa$ is a diffeomorphism of $X$ on a cylinder $Y=\Sigma\times\Lambda$ where
$\Sigma\subset{\Bbb R}^{d-1}$ is an open set and $\Lambda\subset{\Bbb R}$ is an interval.
For  diffeomorphism (\ref{eq:AX1}), function (\ref{eq:IK1}) equals
\begin{equation}
 \tilde{u}(y)= |P_d(x)|^{-1/2} u(x),\quad P_d(x) =\partial P(x)/\partial x_d\neq 0.
\label{eq:AX2}\end{equation}
 We consider $L_2(Y)$ as $L_2(\Lambda;L_2(\Sigma))=:\mbox{\sf H}$. The operator
$\tilde{A}= F_\kappa A F_\kappa^{-1}$ acts in $\mbox{\sf H}$ as multiplication by the independent
variable $\lambda=y_d$. Let us apply Definition~\ref{RS} to the operator
$\tilde{A}$ in the space $\mbox{\sf H}$. Thus we treat $\tilde{A}$ as an
integral operator in the variable $y_d$, and its kernel is an operator acting on functions of the
variable $y_0$. Set ${\cal D}=C_0^\infty(\Sigma)$. If
$\tilde{A}$ has continuous kernel $\tilde{A}^\natural (\mu,\nu):{\cal D}\rightarrow {\cal D}^\prime$
in $\mbox{\sf H}$, then for any $u,v\in C_0^\infty (X)$
\begin{equation}
 (Au,v)_{L_2(X)}= \int_\Lambda \int_\Lambda  (\tilde{A}^\natural (\mu,\nu) \tilde{u}(\nu),
\tilde{v}(\mu))_{L_2(\Sigma)} d\mu d\nu. 
\label{eq:AX3}\end{equation}
This equality makes also sense if
\begin{equation}
X\subset \kappa^{-1}(\Sigma\times \Lambda).
\label{eq:AX3v}\end{equation}
Now, by definition, we accept $\tilde{A}^\natural (\mu,\nu)$ for kernel of the PDO $A$ in the
direct integral associated with the function $P(x)$. This definition depends of course in an obvious
way on choice of the diffeomorphism $\kappa$. 

In the case of PDO,  assumption (\ref{eq:AX3v}) is 
inessential. Indeed, for an arbitrary open set $X$ and $u,v\in C_0^\infty (X)$ consider a partition
of unity $\varphi_n\in C_0^\infty (X)$ such that $\Sigma_{n=1}^N \varphi_n(x)=1$ for $x\in{\rm
supp}\: u\cap {\rm supp}\:v$. Since kernel $k(x,x^\prime)$ (see (\ref{eq:RH1b})) of a PDO
(\ref{eq:RH1}) is a $C^\infty$ function outside of the diagonal,  $\varphi_n A
\varphi_m$ is an integral (in all variables) operator and its kernel is a smooth function provided 
${\rm supp}\: \varphi_n\cap {\rm supp}\:\varphi_m=\emptyset$.  In particular, we note that a PDO $A$
has automatically a smooth kernel $\tilde{A}^\natural (\mu,\nu)$ as long as $\mu\neq \nu$. If 
${\rm supp}\: \varphi_n\cap {\rm supp}\:\varphi_m \neq \emptyset$,  then choosing a sufficiently
fine partition of unity we may achieve that
\[
{\rm supp}\: \varphi_n\cup {\rm supp}\:\varphi_m \subset \kappa^{-1}(\Sigma_{n,m}\times
\Lambda_{n,m})
\]
for diffeomorphism (\ref{eq:AX1}) and suitable choices of a coordinate system  $(x_0,x_d)$ and
sets $\Sigma_{n,m}$, $\Lambda_{n,m}$.

Diagonalization of the operator $P$ requires that $\nabla P(x)\neq 0$ for $x\in X$. Practically,
for the construction of kernel of a PDO $A:C_0^\infty (X)\rightarrow C^\infty (X)$ it suffices to
determine it in a neighbourhood of any  point $x^{(0)}\in X$. Let a unit vector ${\bf n}$ be such
that 
$ <~{\bf n}, \nabla P(x^{(0)}) > \neq 0$. Set $x_d= <x, {\bf n}>$ and let $x_0$ be the orthogonal
projection of $x$ on the hyperplane orthogonal to ${\bf n}$. This defines an admissible coordinate
system $(x_0,x_d)$ in a neighbourhood of the point $x^{(0)}$.

{\bf 5.3.}
In this subsection we set $Y=\Sigma\times\Lambda$ and consider a PDO $\tilde{A}:C_0^\infty
(Y)\rightarrow C^\infty (Y)$ with amplitude $\tilde{\bf a}(y,y^\prime,\eta)$. We shall construct 
kernel
$\tilde{A}(y_d,y_d^\prime):C_0^\infty (\Sigma)\rightarrow C^\infty (\Sigma)$ of $\tilde{A}$ in the
space $L_2(\Lambda; L_2(\Sigma))$ under the assumption that $\tilde{\bf a}$ vanishes in some 
neighbourhood (conical in the variable $\eta$) of the conormal bundle to the hyperplane
$y_d=const$.

More precisely,  assume that there
exists  $\varepsilon>0$ such that
\begin{equation}
  \tilde{{\bf a}}(y,y^\prime,\eta)=0 \quad {\rm if}\quad  |y -y^\prime|\leq \varepsilon \quad {\rm
and}
\quad  |\eta_d|
\geq (1-\varepsilon) |\eta|.
\label{eq:RH2}\end{equation}
Reducing, if necessary, $\Sigma$ and $\Lambda$ we may suppose that (\ref{eq:RH2}) holds for all
$y, y^\prime\in Y$ (not only for $|y -y^\prime|\leq \varepsilon$). Let us check that
$\tilde{A}^\natural(y_d,y^\prime_d):C_0^\infty(\Sigma)\rightarrow C^\infty(\Sigma)$ is the PDO with 
amplitude
\begin{equation}
 \tilde{\bf a}^\natural(y_0,y^\prime_0,\eta_0;y_d,y^\prime_d) =(2\pi)^{-1}\int_{-\infty}^\infty
\tilde{\bf a}(y_0,y_d,y^\prime_0,y^\prime_d,\eta_0,\eta_d)  e^{i<y_d-y^\prime_d,\eta_d>} d\eta_d.
\label{eq:RH6}\end{equation}
Remark, first, that, by assumption (\ref{eq:RH2}), this integral is actually taken over a finite
interval $(-c|\eta_0|, c|\eta_0|)$ where $c=c(\varepsilon) <\infty$. Therefore,
\begin{equation}
\tilde{{\bf a}}^\natural( y_d,y^\prime_d)\in {\cal S}^{m+1}_{\rho,
\delta,\delta}(\Sigma\times \Sigma\times {\Bbb R}^{d-1})
\label{eq:RH3}\end{equation}
 and the corresponding semi-norms are bounded uniformly in $y_d,y^\prime_d$ from any compact
subinterval of $\Lambda$. Moreover, all derivatives
\[ 
(\partial^\alpha_{\eta_0} \partial^{\beta}_{y_0} \partial^{\beta^\prime}_{y^\prime_0} \tilde{{\bf
a}}^\natural)(y_0,y^\prime_0,\eta_0;y_d,y^\prime_d)
\]
 are continuous with respect to $y_d,y^\prime_d$ uniformly in $y_0,y^\prime_0\in K_0$ and $\eta_0,\;
|\eta_0|\leq R_0,$ for any compact $K_0\subset \Sigma$ and any $R_0>0$.

By definition of the PDO  $\tilde{A}^\natural(y_d,y^\prime_d)$, 
\begin{eqnarray*} 
(\tilde{A}^\natural (y_d,y^\prime_d) u(y_d^\prime), v(y_d ))_{L_2(\Sigma)}
\\ =(2\pi)^{-d+1}\int_{\Sigma} \int_{\Sigma} \int_{{\Bbb R}^{d-1}} e^{i<y_0-y_0^\prime,\eta_0>}
 \tilde{{\bf a}}^\natural (y_0,y_0^\prime,\eta_0;y_d,y^\prime_d) u(y_0^\prime,
y_d^\prime)\overline{v (y_0 , y_d )} dy_0 dy_0^\prime d\eta_0,
\end{eqnarray*}
where the amplitude $\tilde{{\bf a}}^\natural$ is given by (\ref{eq:RH6}). Comparing this equality
with (\ref{eq:RH1a})  we arrive at relation  (\ref{eq:RH1c}) for the operator $\tilde{A}$.

Let us summarize the results obtained.

\begin{theorem}\label{RH1}
 Suppose that $\tilde{{\bf a}}\in{\cal S}^{m}_{\rho, \delta,\delta}$ for some $m$ and $\rho>0$,
$\delta<1$ and that condition $(\ref{eq:RH2})$ is satisfied. Then the PDO $\tilde{A}$ with
amplitude  $\tilde{{\bf a}}$ has  continuous kernel
$\tilde{A}^\natural(y_d,y^\prime_d)$ which is also a PDO with amplitude $(\ref{eq:RH6})$ obeying
$(\ref{eq:RH3})$.
 In particular, $\tilde{A}^\natural(y_d,y_d)$ is the PDO with   amplitude
\[
\tilde{\bf a}^\natural(y_0,y^\prime_0,\eta_0;y_d,y_d) =(2\pi)^{-1}\int_{-\infty}^\infty
\tilde{\bf a}(y_0,y_d,y^\prime_0,y_d,\eta_0, \eta_d)   d\eta_d.
\]
\end{theorem}

{\bf 5.4.} 
Let us now consider integral kernels in the direct integral constructed with respect to the operator
of multiplication by some smooth function $P(x)$. Suppose that $\nabla P(x)\neq 0$ in $X$.  For the
construction of kernel of a PDO
$A$, we assume  that its amplitude ${\bf a}(x,x^\prime,\xi)$   equals to zero in some
neighbourhood (conical in the variable $\xi$) of the conormal bundle to  surfaces $S_\lambda$,
that is there exists 
$\epsilon>0$ such that
\begin{equation}
 {\bf a}(x,x^\prime,\xi)=0 \quad {\rm if}\quad  |x -x^\prime|\leq \varepsilon \quad {\rm and} \quad 
|<\xi,\nabla P(x)>| \geq (1-\epsilon) |\xi|\:|\nabla P(x)|.
\label{eq:Rh1}\end{equation}
The variable $x$ in the second condition can of course be replaced by $x^\prime$. Moreover,
reducing, if necessary, $X$  we may suppose that (\ref{eq:Rh1}) holds for all
$x,x^\prime\in X$ (not only for $|x -x^\prime|\leq \varepsilon$).
We may also suppose (see the discussion at the end of subsection 5.2) that, for the diffeomorphism
$\kappa$ defined by (\ref{eq:AX1}), condition (\ref{eq:AX3v}) holds. Then equation of the surface
$S_\lambda=\{ x\in X: P(x)=\lambda\}$ can be written as $x_d=p_\lambda(x_0)$.

Le the diffeomorphism $\kappa$ be defined by equalities (\ref{eq:AX1}). Then
\[
\hskip 1pt ^t \kappa^\prime(x)\eta= \eta_0+   (\nabla  P)(x)\eta_d,\quad \eta=(\eta_0,\eta_d).
\]
For  diffeomorphism (\ref{eq:AX1}), a solution of equation
(\ref{eq:rh2s}) reduces to  construction of a $C^\infty$ vector-function ${\bf q}(x,x^\prime) $ 
satisfying 
\[
  <x-x^\prime, {\bf q}(x,x^\prime)> =P(x) -P(x^\prime).
\]
 Then equation (\ref{eq:rh3}) is fulfilled if 
\begin{equation} 
\xi_0=\eta_0 +{\bf q}_0(x,x^\prime) \eta_d,\quad \xi_d=  q_d(x,x^\prime) \eta_d,\quad {\bf q}=({\bf
q}_0 , q_d).
\label{eq:IK3}\end{equation}
 Necessarily, ${\bf q}(x,x )=(\nabla  P)(x)$ and (cf. (\ref{eq:rh2}))
this function can be constructed by the formula
\[
 {\bf q}(x,x^\prime)=\int_0^1 (\nabla  P)(x+t(x^\prime-x)) dt.
\]
Thus $\tilde{A}$ is the PDO with   amplitude (\ref{eq:rh6}) where $\xi$ is related to $\eta$ by
(\ref{eq:IK3}) and
\[
\alpha (x,x^\prime )=| P_d( x ) P_d( x^\prime ) |^{-1/2} |q_d(x ,x^\prime  )|.
\]

Since $ {\bf q}(x,x)=(\nabla P)(x)$, assumption (\ref{eq:Rh1}) implies that 
  the amplitude ${\bf\tilde{a}}(y, y^\prime,\eta)$    satisfies in $Y=\Sigma\times\Lambda$ 
assumption
(\ref{eq:RH2}). Thus, by Theorem~\ref{RH1},
 the PDO $\tilde{A}$ with this amplitude has  continuous kernel
$\tilde{A}^\natural(y_d,y_d^\prime)$. The operator $\tilde{A}^\natural(y_d,y_d^\prime)$ is  PDO with
the amplitude (\ref{eq:RH6}). So   kernel of the operator $A$ is  also a
PDO with the amplitude  which can be constructed by the formula
\begin{eqnarray}
 {\bf\tilde{a}}^\natural(y_0,y^\prime_0,\eta_0; \mu,\nu) =(2\pi)^{-1}\alpha
(x(\mu), x^\prime(\nu))
\nonumber\\
\times\int_{-\infty}^\infty {\bf a}(x(\mu), x^\prime(\nu), 
\eta_0+ {\bf q}(x(\mu),x^\prime(\nu))\eta_d )  e^{i(\mu-\nu)\eta_d} d\eta_d,
\label{eq:RH6a}\end{eqnarray}
 where 
\[
 x(\mu)= (y_0,p_\mu(y_0)),\quad x^\prime(\nu)= (y_0^\prime,p_\nu(y_0^\prime)).
\]
 Let us formulate the result obtained.

\begin{theorem}\label{AX1}
 Suppose that ${\bf a}\in{\cal S}^{m}_{\rho, \delta,\delta}$ for some $m$ and $\rho>0$, $\delta<1$,
$\rho+\delta\geq 1$ and that condition $(\ref{eq:Rh1})$ is satisfied. Then the PDO $A$ has 
continuous kernel
$\tilde{A}^\natural(\mu,\nu)$ in the direct integral associated with the function $P(x)$. The
operator $\tilde{A}^\natural(\mu,\nu)$  is the PDO with amplitude
$(\ref{eq:RH6a})$ obeying
\[
 \tilde{\bf a}^\natural (\mu,\nu)\in {\cal S}^{m+1}_{\rho, \delta,\delta}(\Sigma\times \Sigma\times
{\Bbb R}^{d-1}).
\]
 In particular, $A^\natural(\lambda,\lambda)$ is the PDO with   amplitude
\begin{equation}
 {\bf\tilde{a}}^\natural(y_0,y^\prime_0,\eta_0; \lambda,\lambda) =(2\pi)^{-1}
\alpha (x(\lambda), x^\prime(\lambda))\int_{-\infty}^\infty {\bf a}(x(\lambda), x^\prime(\lambda), 
\eta_0+ {\bf q}(x(\lambda), x^\prime(\lambda))\eta_d )  d\eta_d.
\label{eq:RH4a}\end{equation}
\end{theorem}

\begin{corollary}\label{AX2}
Let $\psi\in C^\infty ({\Bbb R}) \cap L_1 ({\Bbb R})$, $\int_{-\infty}^\infty \psi(t)dt
=1$ and
\[
\psi_{\varepsilon,\lambda}(x)=\varepsilon^{-1} \psi(\varepsilon^{-1}(P(x)-\lambda)).
\]
 Then for any  $u, v\in C_0^\infty (X)$ the double limit exists 
\begin{equation}
\lim_{\varepsilon,\eta\rightarrow 0} (A \psi_{\varepsilon,\nu}u, \psi_{\eta,\mu}v)_{L_2(X)}=
(\tilde{A}^\natural (\mu,\nu) \tilde{u}(\nu), \tilde{v}(\mu))_{L_2(\Sigma)}
\label{eq:AX5n}\end{equation}
with functions $\tilde{u}(\nu), \tilde{v}(\mu)$ defined by $(\ref{eq:AX2})$.
\end{corollary}

Equality (\ref{eq:AX5n}) allows us to construct the bilinear form of the kernel
 $\tilde{A}^\natural (\mu,\nu)$ in terms of the PDO $A$.

{\bf 5.5.}
Let us consider a particular case $P(x)=x^2$ which is necessary for applications to the scattering
matrix for the Schr\"odinger operator. Define a unitary transformation $W$ of $L_2({\Bbb R}^d)$ on
the space $L_2({\Bbb R}_+;L_2({\Bbb S}^{d-1}))$ of vector-functions $\hat{u}(\lambda;\omega)$ by the
equality
\begin{equation}
\hat{u}(\lambda;\omega)=(Wu)(\lambda;\omega)=2^{-1/2} \lambda^{(d-2)/4} u(\lambda^{1/2}\omega),
\quad \lambda\in {\Bbb R}_+,\; \omega \in {\Bbb S}^{d-1}.
\label{eq:ex1}\end{equation}
Then $WPW^{-1}$ acts as multiplication by the variable $\lambda$ in $L_2({\Bbb R}_+;L_2({\Bbb
S}^{d-1}))$. If the operator $WAW^{-1}$ has continuous kernel $\hat{A}^\natural
(\mu,\nu):C^\infty({\Bbb S}^{d-1})\rightarrow C^\infty({\Bbb S}^{d-1})$, then, by
Definition~\ref{RS},
\begin{equation}
 (Au,v)_{L_2(X)}= \int_0^\infty \int_0^\infty  (\hat{A}^\natural (\mu,\nu) \hat{u}(\nu),
\hat{v}(\mu))_{L_2({\Bbb S}^{d-1})} d\mu d\nu. 
\label{eq:ex2}\end{equation}

To construct kernel $\hat{A}^\natural (\mu,\nu)$ of a PDO $A$ for $\mu,\nu\in (\alpha,\beta)$ we
assume that for some $\varepsilon>0$
\begin{equation}
 {\bf a}(x,x^\prime,\xi)=0 \quad {\rm if}\quad  |x-x^\prime|\leq\varepsilon \quad {\rm and} \quad 
|<\xi,x>| \geq (1-\epsilon)|\xi|\: |x|
\label{eq:ex3}\end{equation}
(and $|x|^2,|x^\prime|^2\in (\alpha,\beta)$). Moreover, we may suppose that ${\bf
a}(x,x^\prime,\xi)=0$ if either $x$ or $x^\prime$ do not belong to a neighbourhood of a point
$x^{(0)}=\lambda^{1/2}_0\omega_0$ where $\lambda _0\in (\alpha,\beta)$, $\omega_0\in{\Bbb
S}^{d-1}$. Therefore condition (\ref{eq:AX3v}) is satisfied for a suitable choice of coordinates
$(x_0,x_d)$ and diffeomorphism  (\ref{eq:AX1}). Clearly, assumption (\ref{eq:ex3}) coincides with
(\ref{eq:Rh1}) for the case $P(x)=x^2$. 
So, according to  Theorem~\ref{AX1}, the PDO $A$ has  continuous kernel
$\tilde{A}^\natural(\mu,\nu)$. The operator $\tilde{A}^\natural(\mu,\nu)$ is also a PDO with
amplitude (\ref{eq:RH6a}) where
\begin{equation}
 {\bf q}(x, x^\prime)=x + x^\prime,\quad \alpha(x,x^\prime)=2^{-1}|x_d x_d^\prime|^{-1/2}\:
|x_d+x_d^\prime|.
\label{eq:ex3n}\end{equation}

This determines kernel $\hat{A}^\natural(\mu,\nu)$. Indeed,
comparing equalities (\ref{eq:AX2}) and (\ref{eq:ex1}), we see that  
$\tilde{u}(\lambda)=Z(\lambda) \hat{u}(\lambda)$, where
\[ (Z(\lambda)f)(y_0)=\lambda^{-(d-1)/2} (\lambda-y_0^2)^{1/4} f( \lambda^{-1/2}y_0, 
(1-\lambda^{-1}y_0^2)^{1/2}).
\]
 Now it follows from  (\ref{eq:AX3})  that
\[ 
\hat{A}^\natural (\mu,\nu)=Z^\ast(\mu) \tilde{A}^\natural (\mu,\nu) Z (\nu)
\] 
satisfies equality  (\ref{eq:ex2}). Thus, under assumption (\ref{eq:ex3}), the PDO $A$ has continuous
kernel $\hat{A}^\natural(\mu,\nu)$. In particular, its diagonal value $\hat{A}^\natural
(\lambda,\lambda)$ is the PDO on ${\Bbb S}^{d-1}$ with the amplitude defined by (\ref{eq:RH4a}) where
${\bf q}$ and $\alpha$ are given by (\ref{eq:ex3n}).


\end{document}